\documentclass[aap,MSNbibl,seceqn,citesort,dvips]{arximspdf}
\usepackage[vtexbibl]{}
\usepackage{graphicx}

\doi{10.1214/10-AAP721}
\volume{21}
\issue{2}
\pubyear{2011}
\firstpage{699}
\lastpage{744}

\makeatletter

\newtheorem{theorem}{Theorem}

\newtheorem{lemma}{Lemma}

\newtheorem{proposition}{Proposition}
\newtheorem{corollary}{Corollary}
\newproclaim{Assumptions}{Assumptions}

\def\eqref#1{(\ref{#1})}
\makeatother

\begin{document}
\begin{frontmatter}

\title{Traveling waves of selective sweeps}
\runtitle{Traveling waves of selective sweeps}

\begin{aug}
\author{\fnms{Rick} \snm{Durrett}\thanksref{t1}\ead
[label=e1]{rtd@math.duke.edu}}
\and
\author{\fnms{John} \snm{Mayberry}\corref{}\thanksref{t2}\ead
[label=e2]{jmayberry@pacific.edu}}

\address{Department of Mathematics \\
Duke University \\
Box 90320 \\
Durham, North Carolina 27708-0320\\ USA \\
\printead{e1}}

\address{Department of Mathematics \\
University of the Pacific \\
3601 Pacific Avenue \\
Stockton, California 95211\\ USA \\
\printead{e2}}

\affiliation{Duke University and University of the Pacific}
\runauthor{R. Durrett and J. Mayberry}

\thankstext{t1}{Supported in part by NSF Grant DMS-07-04996 from the
probability program.}
\thankstext{t2}{Supported in part by NSF RTG Grant DMS-07-39164.}

\end{aug}

\received{\smonth{10} \syear{2009}}
\revised{\smonth{6} \syear{2010}}

\begin{abstract}
The goal of cancer genome sequencing projects is to determine the
genetic alterations that cause common cancers.
Many malignancies arise during the clonal expansion of a benign tumor
which motivates the study of recurrent selective
sweeps in an exponentially growing population. To better understand
this process,
Beerenwinkel et al. [\emph{PLoS Comput. Biol.} \textbf{3} (2007)
2239--2246] consider a Wright--Fisher model in which cells from an exponentially
growing population accumulate advantageous mutations. Simulations show
a traveling wave in which the time of the
first $k$-fold mutant, $T_k$, is approximately linear in $k$ and
heuristics are used to obtain formulas for $E T_k$.
Here, we consider the analogous problem for the Moran model and prove
that as the mutation rate $\mu\to0$, $T_k \sim c_k \log(1/\mu)$,
where the $c_k$ can be computed explicitly. In addition, we derive a
limiting result on a log scale for the size of
$X_k(t) = {}$the number of cells with $k$ mutations at time $t$.
\end{abstract}



\end{frontmatter}

\section{Introduction} \label{intro}

Recent studies have sought to identify the mutations that give rise to
common cancers by sequencing protein-coding genes in common tumor
types, including breast and colon cancer \cite{SJW,WPJ},
pancreatic cancer \cite{Jetal} and glioblastoma \cite{Petal,TCGA}. The last study is part of a 100 million dollar pilot
project of the NIH, which could lead to a 1.5 billion dollar effort.
These studies have rediscovered genes known to play a role in cancer
(e.g., APC, KRAS and TP53 in colon cancer), but they have also found
that tumors contain a large number of mutations. Analysis of 13,023
genes in 11 breast and 11 colorectal cancers in Sjoblom et al.~\cite
{SJW} revealed that individual tumors accumulate an average of $\approx
$90 mutated genes, but only a subset of these contribute to the
development of cancer.

Follow-up work in Wood et al. \cite{WPJ} studied 18,191 distinct genes
in the same 22 samples. Any gene that was mutated in a tumor but not
normal tissue was analyzed in 24 additional tumors, and selected genes
were further analyzed in 96 colorectal cancers. Statistical analysis
suggested that most of the $\approx$80 mutations in an individual
tumor were harmless and that $<$15 were likely to be responsible for
driving the initiation, progression or maintenance of the tumor. These
two types of mutations are commonly referred to as ``drivers'' and
``passengers.'' The latter provide no selective advantage to the growing
cancer mass, but are retained by chance during repeated rounds of cell
division and clonal expansion (exponential growth).

The results of \cite{SJW} and \cite{WPJ} are in contrast with the
long-held belief that most cancers are the end result of a handful of
mutations. Armitage and Doll \cite{AD} constructed log--log plots of
cancer mortality versus age and found slopes of 5.18 and 4.97 for colon
cancer in men and women, respectively. From this they predicted that
the occurrence of colon cancer was the result of a six-stage process.
In essence their argument is that the density function of the sum of
six exponentials with rates $\mu_i$ is
\[
\approx\mu_1 \cdots\mu_6 t^5/5! \qquad\mbox{for small $t$.}
\]
This result yields the density of the well-known gamma distribution
when all the $\mu_i$ are equal, but only readers with well developed
skills in calculus (or complex variables) will succeed in deriving this
result for unequal $\mu_i$.

The work in \cite{AD} and the subsequent work of Knudson \cite{K0}, who used
statistics to argue that retinoblastoma was the end result of two
mutations, gave rise to a large amount of work; see  \cite
{KNRC} and the books by Wodarz and Komarova \cite{WK} and Frank \cite
{Fr} for surveys. From this large body of work on multistage
carcinogenesis, we will cite only two sources. Luebeck and Moolgavakar
\cite{LM} used multistage models to fit the age-specific incidence of
colorectal cancers in the SEER registry (which covers 10\% of the
US~population) to conclude that a four-stage model gives the best
fit. Calabrese et al.~\cite{Cetal} used data for 1022 colorectal
cancers to argue that ``sporadic'' cancers developed after six
mutations, but that in the subgroup of individuals with strong familial
predispositions, only five mutations were required.

There is good reason to doubt some of the conclusions of \cite{SJW}
and \cite{WPJ}. First, the statistical methods of \cite{SJW} have
been criticized (see letters on pages 762--763 in the February 9, 2007
issue of \emph{Science}). Furthermore, in \cite{WPJ}, a follow-up
study on 40 of the 119 highest scoring genes, chosen because they were
in pathways of biological interest, showed that 15 of the 40 genes
(37.5\%) were not mutated in any of the 96 tumors, casting doubt on the
claimed 10\% false discovery rate. However, the more recent studies
\cite{Jetal,Petal,TCGA} using well-known and
trusted statistical methods have found similar patterns: an average of
63 mutations in pancreatic cancers and 47 in glioblastoma.

To better understand this process by which an exponentially growing
cell mass accumulates driver and passenger mutations, and, in
particular, to understand the data in \cite{SJW}, Beerenwinkel et
al.~\cite{BN} considered a Wright--Fisher model with selection and
mutation in an exponentially growing population. They assumed that
there were 100 potential driver genes and asked for the waiting time
until one cell has accumulated $k$ mutations. Simulations (see their
Figure 3) showed that a traveling wave developed, in which the time
until the first $k$-fold mutant was approximately linear in $k,$ and
the authors used heuristic arguments to obtain quantitative predictions
for the first time that a cell with $k$ mutations appears.

Here we will consider this problem for the analogous Moran model and
prove asymptotic results as the mutation rate $\mu\to0$ for the
behavior of $X_j^\mu(t) = {}$the number of cells with $j$ mutations at
time $t$. A cell with $j$ mutations will be referred to as a type-$j$
individual. Our main result is Theorem \ref{induct-grow}, which allows
for an exponentially growing population $N^\mu(t)$ of individuals. The
process of fixation of advantageous mutations in a population of
constant size has been the subject of much theoretical work (see, e.g.,
Chapter 6 of  \cite{RD}), so it is natural to ask how the
behavior changes in an exponentially growing population. A second
difference from the standard theory of the fixation of a single
mutation is that we consider a situation in which new mutations arise
before older ones have gone to fixation, a~process often referred to as
``stochastic tunneling.'' The resulting ``Hill--Robertson'' interference
(see, e.g., Section 8.3 in \cite{RD}) can be analyzed here because
only the newest mutation is stochastic while the older mutations behave
deterministically. This idea was used by Rouzine et al.~in \cite{RWC}
(and later developed in more detail in \cite{BRW,RBW}) as a
heuristic, but here it leads to rigorous results.

The rest of the paper is organized as follows. In Section \ref
{sec-fixedpop} we begin with a fixed population size of $N = \mu
^{-\alpha}$ individuals and state Theorem \ref{thm-mainlimit}, which
says that when time is scaled by $L=\log(1/\mu)$, the $\log$ sizes
of $X_j$, divided by $L$, converge to a limit that is deterministic and
piecewise linear and so the time the first type-$j$ individual appears
is $O(\log(1/\mu))$. Since we have assumed the population size is
$\mu^{-\alpha}$, this time scale agrees with (i) results in
\cite{YE,YEC}, which show that the rate of adaptation
(defined as the change in the mean fitness of the population) for a
related fixed-population-size Moran model is $O(\log N)$ and (ii)
simulations in  \cite{DF} which suggest that the speed
of adaptation depends logarithmically on both the mutation rate and the
population size. In Section \ref{sec-grow} we return to the growing
population scenario and state our main result, Theorem \ref
{induct-grow}, which generalizes Theorem~\ref{thm-mainlimit}. Section
\ref{examples} contains some examples elucidating the nature of the
limit in Theorem~\ref{induct-grow} and illustrating the
traveling-wave-like behavior of the limit. We state the main tools used
to prove Theorem~\ref{induct-grow} in Section \ref{pf-ideas}, and
Sections \ref{sec-startup} and \ref{sec-inductpf} contain the technical details.

\subsection{Fixed population size: Main result} \label{sec-fixedpop}

We begin by considering our Moran model in a fixed population of $N$
individuals and return to our analysis of the exponentially growing
population in Section \ref{sec-grow}. We assume that:

\begin{longlist}[(iii)]
\item[(i)] initially, all individuals are of type $0$;

\item[(ii)] type-$j$ individuals mutate to individuals of type $j+1$ at rate
$\mu$;

\item[(iii)] all individuals die at rate 1 and, upon death, are replaced by an
individual of type $j$ with probability
\[
\frac{(1+\gamma)^j X_j^\mu(t)}{W^\mu(t)},
\]
where $(1+\gamma)^j$ is the relative fitness of type-$j$ individuals
compared to type-0, and
\begin{eqnarray*}
W^\mu(t) = \sum_{i=0}^{\infty} (1+\gamma)^i X_i^\mu(t)
\end{eqnarray*}
is the ``total fitness'' of the population. We assume throughout that
$\gamma>0$ is fixed (i.e., mutations are advantageous). Approximations of
the time the first type-$k$ individual appears have been carried out
for the neutral case ($\gamma= 0$) in \cite{INM,HIM,DSS,JS}
(and applied to regulatory sequence evolution in
\cite{DS}). The case $\gamma<0$ is of interest in studying Muller's
ratchet \cite{Mul}, but since deleterious mutations behave
very differently from advantageous mutations, we will not consider this
case here.
\end{longlist}

We will suppose throughout that $N \gg1/\mu$, that is, $N\mu\to
\infty$. If $N\mu\to0$, then the 1's arise
and go almost to fixation before the first mutation to a 2 occurs, so
the times between fixations are independent exponentials. We will not
here consider the borderline scenario, although we note that in the
case when $N \mu\to c_1 >0$ and $N \gamma\to c_2 >0$, the limiting
behavior of the system has been well studied and we obtain a diffusion
limit describing the evolution of type-$j$ frequencies (see, e.g.,
Sections 7.2 and 8.1 in  \cite{RD}). Let $T_0^\mu=0,$ and for
$j\ge1$ define
\begin{eqnarray*}
T_j^\mu&= \inf\{t \geq0\dvtx X_j^\mu(t) \geq1\},
\end{eqnarray*}
that is, $T_j^\mu$ is the time of the first appearance of a type-$j$
individual. In order to study the birth times $T^\mu_j$ we will prove
a limit theorem for the sizes of the $X_j^\mu(t)$ on a $\log$ scale.
Let $\log^+x = \max\{\log x, 0 \}$, $L=\log(1/\mu)$ and define
\[
\gamma_j = (1+\gamma)^j -1
\]
for all $j \in\mathbb{Z}$. In what follows, we shall use $\lfloor
x\rfloor$ to denote the greatest integer less than or equal to $x$ and
let $\lceil x \rceil= \lfloor x\rfloor+ 1$.

\begin{theorem} \label{thm-mainlimit}
Suppose that $X_0^\mu(0)=N$ with $N = \lceil\mu^{-\alpha} \rceil$
for some $\alpha>1$ and suppose that $\gamma\in G(\alpha)$, the set of
generic parameter values defined in \eqref{generic}. Then, as $\mu\to0,$
\begin{eqnarray*}
Y^\mu_j(t) \equiv\frac{1}{L} \log^+\bigl(X_j(Lt/\gamma)\bigr) \to y_j(t)
\qquad\mbox{in probability}
\end{eqnarray*}
uniformly on compact subsets of $(0,t^*),$ where $t^* = t^*(\alpha
,\gamma
) > 0 $ is defined in \eqref{tstar}. The limit $y_j(t)$, which depends
on the parameters $\alpha,\gamma$, is deterministic and piecewise linear
and will be described by \textup{(a)} and \textup{(b)} below. Furthermore, if we define
\[
t_j = t_j(\alpha,\gamma) = \inf\{t\dvtx y_j(t) = 1\}
\]
for $j \geq0$, then
\[
\frac{T_{j+1}^\mu}{L/\gamma} \to t_j \qquad\mbox{in probability}
\]
as $\mu\to0$ for all $j \geq0$.
\end{theorem}

\begin{enumerate}[(a)]
\item[(a)] {\it Initial behavior.} $y_j(0)= (\alpha-j)^+$. The convergence
only occurs on $(0,\infty)$
because we have $Y_j^\mu(0) =0 $ for all $j \geq1$ by assumption, so
a discontinuity is created at time 0.

\item[(b)] {\it Inductive step.} Let $s_0=0$ and suppose that at some time
$s_n \geq0,$ the following conditions are satisfied:
\begin{longlist}[(ii)]
\item[(i)] $m_n \equiv\max\{j\dvtx y_j(s_n)=\alpha\}$ and $k_n \equiv
\max\{j\dvtx y_j(s_n) >0\}$ both exist and $y_j(s_n) >0$ for all $m_n \leq
j \leq k_n;$
\item[(ii)] $y_{j+1}(s_{n}) \geq y_j(s_n) -1$ for all $0 \leq j \leq
k_n$ so that, in particular, \mbox{$y_{k_n}(s_n) \leq1$.}
\end{longlist}
\end{enumerate}
Let $k_n^*=k_n$ if $y_{k_n}(s_n)<1$, $k_n^*=k_n+1$ if $y_{k_n}(s_n)=1$
and define\vspace*{-2pt}
\[
\delta_{n,j} =
\cases{ \bigl(\alpha-y_j(s_n)\bigr)\gamma/\gamma_{j-m_n}, &\quad $m_n < j <
k_n^*$,
\cr
\bigl(1-y_{k_n^*}(s_n)\bigr)\gamma/\gamma_{k_n^*-m_n}, &\quad $j=k_n^*$.
}\vspace*{-2pt}
\]
For all $t \le\Delta_n \equiv\min\{\delta_{n,j}\dvtx m_n < j \leq
k_n^*\}$, we then have\vspace*{-2pt}
\begin{eqnarray*}
y_j(s_n + t) =
\cases{ \bigl(y_j(s_n) + t\gamma_{j-m_n}/\gamma\bigr)^+, &\quad $j \leq k_n^*$,
\cr
0, &\quad $j > k_n^*$,
}\vspace*{-2pt}
\end{eqnarray*}
and conditions (i) and (ii) are satisfied at time $s_{n+1} = s_n +
\Delta_n$.

Our description of the limiting dynamical system can be understood as
follows. If type $m_n$ is the most fit of the dominant types in the
population at time $s_n$, then the $y_j(s_n +t)$, $m_n < j \leq k_n^*$,
grow linearly with slope $\gamma_{j-m_n}/\gamma> 0,$ while the
$y_j(s_n+t)$, $j<m_n$, decrease linearly with slope $\gamma
_{j-m_n}/\gamma< 0$ until they hit zero and $y_{m_n}(s_n+t)$ stays
constant. These rates remain valid until either $y_j$ reaches level
$\alpha$ for some $m_n < j < k_n^*$ and there is a change in the most
fit dominant type, or $y_{k_n^*}$ reaches level 1 and individuals of
type $k_n^*+1$ are born. These two events correspond to $\Delta_n =
\delta_{n,j}$ and $\Delta_n = \delta_{n,k_n^*}$, respectively. The
condition $y_{j+1}(s_n) \ge y_j(s_n) - 1$ guarantees that after birth,
the growth of type-$j$ individuals is driven by selection and not by
mutations from type-$(j-1)$ individuals. If this condition failed, we
would encounter a discontinuity in the limiting dynamics like the one
at time 0. We have rescaled time by $\gamma^{-1}$ since, in most cases of
interest, $\gamma$ is small (e.g., $\gamma<0.01$) and when $\gamma$ is
small, we have $\gamma_j/\gamma\approx j$ so that the limit process described
above is almost independent of $\gamma$.

Parts (a) and (b) together describe the limiting dynamical system for all times\vspace*{-2pt}
%
\begin{equation}\label{tstar}
t < t^* \equiv\sum_{n=1}^\infty\Delta_n\vspace*{-2pt}
\end{equation}
since by part (a), (i) and (ii) in (b) hold at time $s_0=0,$ and it is
easy to see from the form of $y_j(t)$ that if (i) and (ii) hold at time
$s_n$, then they also hold at time $s_{n+1}=s_n+\Delta_n$. Note that
the form of the limit implies that at times $t \notin\{s_n\}_{n \geq
0}$, there is always a unique $j$ such that $y_j(t) = \alpha$, that
is, a unique dominant type. This observation will prove useful on a
number of occasions. In Section \ref{examples} we shall see examples
in which $t^* = \infty$, but in Section \ref{sec-bu} we will prove
the following result.\vadjust{\goodbreak}

\begin{lemma} \label{blowup}
For any $\gamma>0$, there exists $\alpha_\gamma$ such that $t^* <
\infty$
whenever $\alpha> \alpha_\gamma$.
\end{lemma}

 Therefore, in general, our construction cannot be extended up to
arbitrarily large times. We prove this lemma by showing that if $\alpha
$ is large, then an infinite number of types will be born before any
existing type-$j \geq1$ individual achieves dominance. However, since
it is easy to see by construction that we have $t^* \geq t_j$ for any
$j$, this is the only way that blow-up can occur. Hence, the limit
process will accumulate an infinite number of mutations before time
$t^*,$ which means our approximation is valid on any time interval of
practical importance.


The generic set of parameter values is
%
\begin{equation}\label{generic}
G(\alpha) \equiv\{\gamma\in(0,\infty)\dvtx \delta_{n,j} \neq\delta
_{n,i} \mbox{ for all } i \neq j, n \geq0 \}.
\end{equation}
In other words, when $\gamma\in G(\alpha)$, there is always a unique
$j_n$ such that $\Delta_n = \delta_{n,j_n}$. For any given value of
$\alpha$, $G(\alpha)$ is clearly countable, so our result is good
enough for applications. If $\delta_{n,i} = \delta_{n,j}$ for some $i
\neq j$, then either (i) type $i$'s and type $j$'s reach level $\alpha
$ at the same time or (ii) type $i$'s reach level $\alpha$ at the same
time that type $j$'s reach level 1. It is tempting to argue that since
generic parameters are dense, the result for general parameters
follows, but proving this is made difficult by the fact that the growth
rates are not continuous functions of the parameters since they depend
on $m_n = \max\{ j \dvtx y_j(s_n)=\alpha\}$.

\begin{figure}

\includegraphics{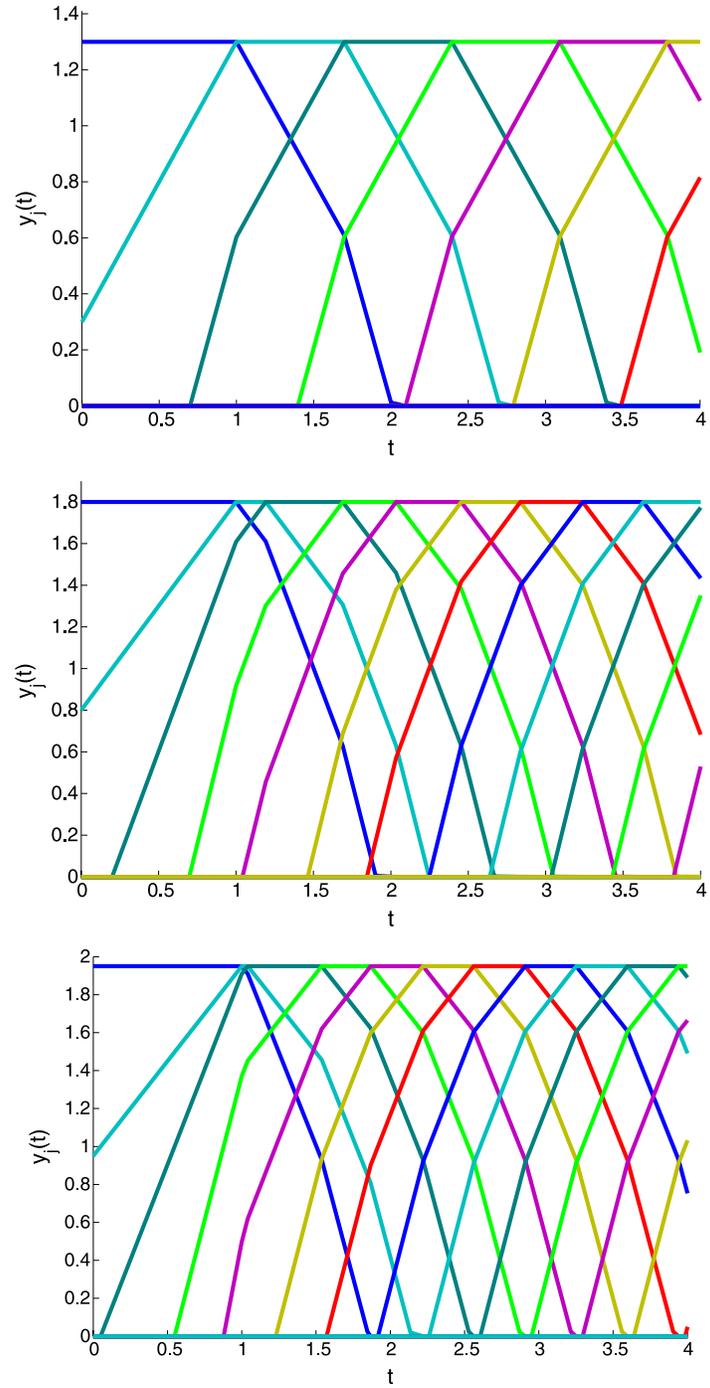}

  \caption{Plot of the limiting dynamical system in Theorem \protect\ref{thm-mainlimit} (fixed population size).
  Parameters: $\gamma=0.01$, $\alpha = 1.3, 1.8, 1.95$ (top to bottom).}\label{fig:r123}
\end{figure}

Theorem \ref{thm-mainlimit} is very general, but not very transparent.
In Section \ref{examples} we will give some examples in which more
explicit expressions for the birth times $t_j$ are available. Figure
\ref{fig:r123} shows examples in the first three ``regimes'' of
behavior that we will consider. In the $j$th regime, type $m+j$ arises
(but not type $m+j+1$) before type~$m$ ``fixates,'' that is, is of
order $N = \lceil\mu^{-\alpha}\rceil$. These regimes closely
correspond to the different scenarios considered in \cite
{BRW}, in which the ``stochastic edge,'' that is, the class of the most
fit mutant present at positive quantities, is always assumed to be $q$
fitness classes ahead of the population bulk. $q$ is referred to as the
``lead.'' In the notation of Theorem \ref{thm-mainlimit}, the lead is
always $k_n^*$ on the interval $[s_n,s_{n+1}),$ and in regime $j$, the
lead is always $j$. In all three examples in Figure \ref{fig:r123}, we
see the traveling-wave-like behavior observed in the simulations of
Beerenwinkel et al.~\cite{BN} (see also \cite{RWC}).
The time between successive waves is constant in the example from
regime 1, while in the examples from regimes 2 and 3, the time between
successive waves is not constant, but converges to a constant as the
number of waves goes to infinity.

\subsection{Growing population} \label{sec-grow}

We now consider a growing population of individuals $N^\mu(t)$, $t
\geq0,$ with a random initial population in $\mathbb{N} = \{1,2,\ldots\}
$ distributed according to some measure $\nu_0$. At time $0$, all
individuals are of type 0 and we suppose that, in addition to the
previously imposed Moran dynamics, at rate $\rho N^\mu(t)$, $\rho\geq
0$, new individuals are added and their type is chosen to be $k$ with
probability
\[
\frac{(1+\gamma)^k X_k^\mu(t)}{W^\mu(t)}.
\]
As in the case of fixed population size, we are able to derive a
limiting, piecewise linear approximation to
\[
Y_k^\mu(t) \equiv(1/L) \log^+ X_k^\mu(Lt/\gamma).
\]
To determine the correct growth rates, suppose that there are $x_j$
individuals of type~$j$ and that the population size is $N$. We then
have the jump rates
\begin{eqnarray*}
x_j &\mapsto& x_j +1 \qquad\mbox{rate: } [(1+\rho) N - x_j]
\frac{(1+\gamma)^jx_j}{\sum_{i \geq0} (1+\gamma)^i x_i} + \mu
x_{j-1},
\\
x_j &\mapsto& x_j -1 \qquad\mbox{rate: }  x_j \frac{ \sum_{i
\neq j} (1+\gamma)^i x_i }{\sum_{i \geq0} (1+\gamma)^i x_i } -\mu x_j.
\end{eqnarray*}
If mutations can be ignored, then the growth rate of type $j$'s is
\begin{eqnarray*}
\frac{\sum_{i \geq0} [(1+\rho)(1+\gamma)^j-(1+\gamma)^i]x_ix_j
}{\sum_{i
\geq0} (1+\gamma)^i x_i} \approx[(1+\rho)(1+\gamma)^{j-m}-1] x_j
\end{eqnarray*}
if $x_i = o(N)$ for $i\neq m$ (recall that in the limiting dynamical
system from Theorem~\ref{thm-mainlimit}, there is a unique dominant
type at time $t$ for all but a countable number of times). This yields
the expression
\[
\lambda_{j-m} \equiv(1+\rho)(1+\gamma)^{j-m}-1
\]
for the limiting growth rate of type $j$'s in a population dominated by
type $m$.

If type-$j$ individuals have size $(1/\mu)^x$ at time 0, are growing
at rate $\lambda_{\ell(j)}$ for some $\ell(j) \geq1$ and the
initial total population size is $(1/\mu)^z$, then type $j$'s will
achieve fixation at the approximate time $t$ satisfying
\[
(1/\mu)^x e^{\lambda_{\ell(j)} t} = (1/\mu)^z e^{\rho t}
\quad\mbox{or}\quad t = \frac{z-x}{\lambda_{\ell(j)}-\rho} \log
(1/\mu).
\]
This leads to the following result. Theorem \ref{thm-mainlimit}
corresponds to the special case $\rho=0$.

\begin{theorem} \label{induct-grow}
Let $F^\mu(t) = (1/L) \log N^\mu(tL/\gamma)$ and suppose that $F^\mu(0)
\to\alpha$ in probability for some $\alpha> 0$. Then, for all
$\gamma
\in G(\alpha,\rho)$, the set of generic parameter values given in
\eqref{genericrho}, as $\mu\to0$ we have $F^\mu(t) \to\alpha+ t
\rho/\gamma$ and $Y^\mu_j(t) \to y_j(t)$ in probability uniformly on
compact subsets of $[0,\infty)$ and $(0,t^*),$ respectively, where
$t^* = t^*(\alpha,\rho,\gamma)$ is given in \eqref{tstarrho}. The
limits $y_j(t)$, which depend on the parameters $(\alpha,\rho,\gamma)$,
are deterministic and piecewise linear and described by \textup{(a)} and \textup{(b)}
below. Furthermore, if we define
\[
t_j = t_j(\alpha,\rho,\gamma) = \inf\{t\dvtx y_j(t) = 1\}
\]
for $j \geq0$, then
\[
\frac{T_{j+1}^\mu}{L/\gamma} \to t_j \qquad\mbox{in probability}
\]
as $\mu\to0$ for all $j \geq0$.
\end{theorem}

\begin{enumerate}[(a)]
\item[(a)] {\it Initial behavior.} $y_j(0)= (\alpha-j)^+$.

\item[(b)] {\it Inductive step.} Let $s_0=0$ and suppose that at some time
$s_n \geq0$ the following conditions are satisfied:
\begin{longlist}[(ii)]
\item[(i)] $m_n \equiv\max\{j\dvtx y_j(s_n)=\alpha+ \rho s_n \}$ and
$k_n \equiv\max\{j\dvtx y_j(s_n) >0\}$ both exist and $y_j(s_n) >0$ for
all $m_n \leq j \leq k_n;$
\item[(ii)] $y_{j+1}(s_{n}) \geq y_j(s_n) -1$ for all $0 \leq j \leq
k_n$ so that, in particular, \mbox{$y_{k_n}(s_n) \leq1$.}
\end{longlist}
\end{enumerate}
Let $\alpha_n = \alpha+ \rho s_n$, $k_n^*=k_n$ if $y_{k_n}(s_n)<1$,
$k_n^*=k_n+1$ if $y_{k_n}(s_n)=1$ and define\vspace*{-2pt}
\[
\delta_{n,j} =
\cases{ \bigl(\alpha_n-y_j(s_n)\bigr)\gamma/(\lambda_{j-m_n}-\rho), &
\quad $m_n < j < k_n^*$,
\cr
\bigl(1-y_{k_n^*}(s_n)\bigr)\gamma/\lambda_{k_n^*-m_n}, &\quad $j=k_n^*$.
}\vspace*{-2pt}
\]
For all $t \le\Delta_n \equiv\min\{\delta_{n,j}\dvtx m_n < j \leq
k_n^*\}$, we then have\vspace*{-2pt}
\begin{eqnarray*}
y_j(s_n + t) =
\cases{ \bigl(y_j(s_n) + t\lambda_{j-m_n}/\gamma\bigr)^+, &\quad $j \leq k_n^*$,
\cr
0, &\quad $j > k_n^*$,
}\vspace*{-2pt}
\end{eqnarray*}
and conditions (i) and (ii) are satisfied at time $s_{n+1} = s_n +
\Delta_n$.

The generic set of parameter values is\vspace*{-2pt}
%
\begin{equation}\label{genericrho}
G(\alpha,\rho) \equiv\{\gamma\in(0,\infty)\dvtx \delta_{n,j} \neq
\delta
_{n,i} \mbox{ for all } i \neq j, n \geq0 \}\vspace*{-2pt}
\end{equation}
and, of course,\vspace*{-2pt}
%
\begin{equation}\label{tstarrho}
t^* \equiv\sum_{n=1}^\infty\Delta_n.\vspace*{-2pt}
\end{equation}
The argument which we use to prove Lemma \ref{blowup} implies that
$t^* < \infty$ whenever $\rho>0$ since $\alpha_n \to\infty$ as $n
\to\infty$.

\begin{figure}

\includegraphics{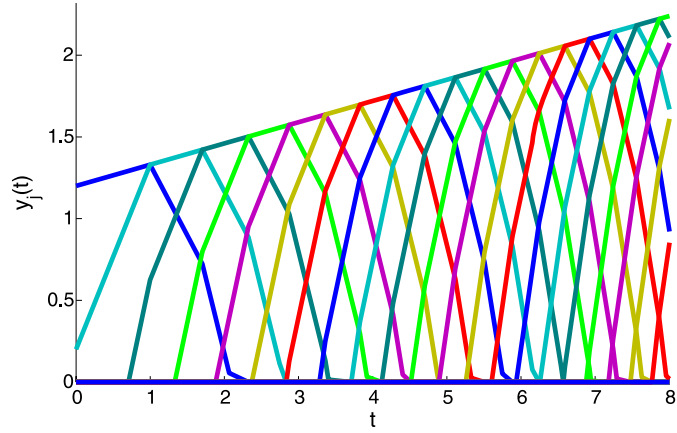}

  \caption{Plot of the limiting dynamical system in Theorem \protect\ref{induct-grow}
  (growing population size). Parameters: $\rho = 0.0013$, $\gamma = 0.01$ and $\alpha=1.2$.}
  \label{fig:grow}
\end{figure}

An example is given in Figure \ref{fig:grow}. Since the population
size is growing, we progress through the different ``regimes'' of
behavior defined earlier for the fixed population size, and the time
between successive waves of sweeps decreases. This behavior can also be
seen in Figure 3 of \cite{BN}. Here we are dealing
with the small mutation limit so that our waves have sharp peaks.

Motivated by the statistical analysis of cancer data in \cite{SJW},
Beerenwinkel et al.~\cite{BN} were interested in the time $T_{20}^\mu
$ at which a cell first accumulates 20 mutations. Their choice of the
number 20 was inspired by data from \cite{SJW}. Using heuristics they
obtained the approximation\vspace*{-2pt}
%
\begin{equation}\label{Bapprox}
T_j^\mu\approx s_j = j \frac{(\log(\gamma/\mu))^2}{\gamma\log
(N(0)N(T_{20}^\mu))}\vspace*{-2pt}
\end{equation}
for $j \leq20$. Their model evolves in discrete time, but the
heuristics only use the fact that the drift in the Wright--Fisher
diffusion limit (ignoring mutations) is given by\vspace*{-2pt}
\[
b_j(x) \approx\gamma x_j(j - \langle j\rangle),
\]
\begin{figure}[b] 

\includegraphics{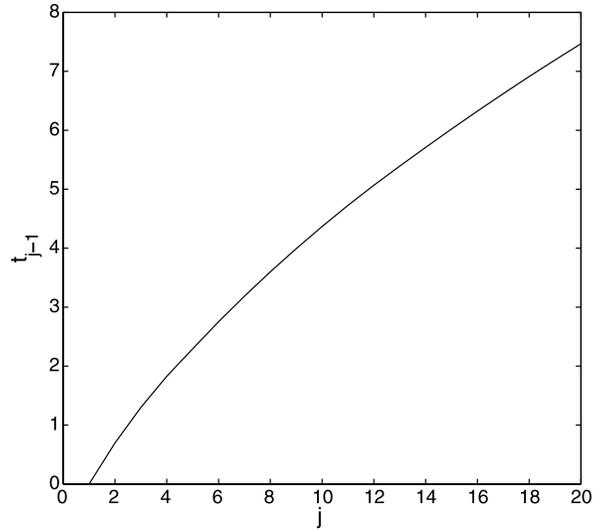}

  \caption{Plot of $t_{j-1}$, the constants from Theorem \protect\ref{induct-grow}, as a
  function of $j$. Same parameter values as in Figure \protect\ref{fig:grow}.}
  \label{fig:jvsSj}
\end{figure}
where $\langle j\rangle = \sum j x_j$ (see  \cite{RD}, page 253). To get the
same drift in continuous time, we need to rescale time by $2/N$, as
opposed to $1/N,$ and hence we should replace $\gamma$ by $\gamma/2$
and $\mu
$ by $\mu/2$ to obtain the analogous approximations for the Moran
model. The important point to note is that the approximation in \eqref
{Bapprox} is linear in $j$ and hence yields constant estimates for the
increments $T_{j}^\mu-T_{j-1}^\mu$, whereas we can see that in the
limiting dynamical system, the increments are not constant, but
decrease in length as the population size increases. Figure \ref
{fig:jvsSj} shows a plot of $t_{j-1}$, the constants from Theorem \ref
{induct-grow}, as a function of $j$ and illustrates the nonlinearity in $j$.


\section{Examples} \label{examples}

In this section, we discuss some examples with constant population size
($\rho= 0$) in which the explicit computation of the times $t_j$ is
possible. To do so, it is more convenient to study the increments
\begin{eqnarray*}
\tau_j^\mu\equiv T_j^\mu- T_{j-1}^\mu.
\end{eqnarray*}
Theorem \ref{thm-mainlimit} then implies that
\[
\frac{\tau_j^\mu}{L/\gamma} \to\beta_j,
\]
where $\beta_j = t_{j-1}-t_{j-2}$ for all $j \geq1$ if we use the
convention that $t_{-1}=0$. We begin with the first regime of behavior
where the lead is always 1 and we have $\beta_j \equiv\beta$ for all
$j \geq2$. In Section \ref{sec-r2} we move on to discuss regime 2
where the lead is always 2. In this case, we will see that $\beta_j$
depends on $j$, but we have $\beta_j \to\beta$ as $j \to\infty$ so
that asymptotically, the times between successive waves is constant.
Section \ref{sec-r3} contains some conjectures on the parameter range
for the regime $j \geq3$. In all these scenarios, we conjecture that
$t^* = \infty$ and our limiting result holds for all time. In Section
\ref{sec-bu}, we will prove Lemma \ref{blowup}, showing that for any
$\gamma>0$, there exists $\alpha_\gamma$ such that $t^* < \infty$ whenever
$\alpha> \alpha_\gamma$.

\subsection{Results for regime 1} \label{sec-r1}

Let $r_2 = 1 + \gamma/\gamma_2$. The first regime occurs for $1<
\alpha< r_2$. If $\gamma$ is small, then $\gamma/\gamma_2 \approx
1/2$ and the condition is roughly $\alpha\in(1,3/2)$. If $\gamma>0$,
then $\gamma/\gamma_2 = 1/(2+\gamma) < 1/2,$ so $\alpha< 2$
throughout regime 1.

\begin{table}[b]
\tablewidth=295pt
\caption{Sizes in regime 1. Times are given in units of $L/\gamma$,
entries are the size given as a power of $1/\mu$ and 0 indicates when
the first of the type is born. The first row comes from \textup{(a)}, the next
four from applications of \textup{(b)}}
\label{table:r1}
\begin{tabular*}{295pt}{@{\extracolsep{4in minus 4in}}lccccc@{}}
\hline
\textbf{Time} & \textbf{Time increment} & \textbf{Type 0} & \textbf{Type 1} & \textbf{Type 2} & \textbf{Type 3} \\
\hline
$0+$ & & $\alpha$ & $\alpha-1$ \\
$s_1$ & $\Delta_0 = 2-\alpha$ & $\alpha$ & 1 & 0 \\
$s_2$ & $\Delta_1 = \alpha-1$ & $\alpha$ & $\alpha$ & $\gamma_2
(\alpha-1)/\gamma$ \\
$s_3$ & $\Delta_2 = 1 - \Delta_1\gamma_2/\gamma$ & & $\alpha$ & 1
& 0\\
$s_4$ & $\Delta_3 = \alpha-1$ & & $\alpha$ & $\alpha$ & $\gamma_2
(\alpha-1)/\gamma$\\
\hline
\end{tabular*}
\end{table}

Table \ref{table:r1} summarizes the situation. To explain the entries,
we note that applying part (a) of the limit description implies that
$y_1(0)=\alpha-1,$ and part (b) then implies that
\[
y_1(s) = (\alpha- 1) + s
\]
for $s \leq\Delta_0 = 2-\alpha$. Since we have assumed that $\alpha
< r_2$, we have
\[
\Delta_1 = \delta_{1,1} \wedge\delta_{1,2} = (\alpha-1) \wedge
\frac{\gamma_2}{\gamma} = \alpha-1,
\]
and applying part (b) tells us that we have $y_2(s_1+t) = t\gamma
_2/\gamma$
for all $t \leq\Delta_1$. Another application of (b) then yields
$\Delta_2 = \delta_{2,2} = 1 - y_2(s_2),$ which gives the additional
amount of time needed for $y_2$ to hit $1$. Since the relative sizes of
1's, 2's and 3's at time $s_3$ are the same as the relative sizes of
0's, 1's and 2's at time $s_1$, we obtain the following result giving
the limiting coefficients of $\tau_j^\mu$.

\begin{corollary} \label{cor-twos}
Suppose that $N=\mu^{-\alpha}$ for some $\alpha\in(1,r_2)$. Then,
as $\mu\to0,$
\begin{eqnarray*}
\frac{\tau_1^\mu}{L/\gamma} \to(2-\alpha), \quad\mbox{and for all
$j \geq2,$}\quad
\frac{\tau_j^\mu}{L/\gamma} \to\beta\qquad\mbox{in probability,}
\end{eqnarray*}
where $
\beta\equiv\Delta_1 + \Delta_2 = (\alpha-1)+ 1-(\alpha-1)\frac
{\gamma_2}{\gamma} = (2+\gamma) - (1+\gamma)\alpha.
$
\end{corollary}

Figure \ref{fig:r123} illustrates the limiting dynamical system in the
case where $\gamma=0.01$ and $\alpha= 1.3$. We can see that in regime 1,
the system is characterized by a ``traveling wave of selective sweeps''
in type space, that is, the growth and decay of types $j \ge2$ occur
translated in time by a fixed amount. In Figure \ref{fig:typedist} we
show the distributions of types at the times when type-5, -9, -13 and
-17 individuals are born (the times $t_4$, $t_8$, $t_{12}$ and $t_{16}$
from Theorem \ref{thm-mainlimit}). As we move from time $t_j$ to
$t_{j+4}$, the distribution is shifted by a constant amount.

\subsection{Results for regime 2} \label{sec-r2}

Regime 2 occurs for $r_2 < \alpha< r_3$ with $r_3 = r_2 + \gamma
/\gamma_3$. When $\gamma$ is small, $\gamma/\gamma_3 \approx1/3,$
so this regime is roughly $\alpha\in(3/2,11/6)$. In general, $r_3 <
11/6,$ so we have $\alpha< 2$ throughout this regime. As in the
previous section, it is easiest to explain the conclusions of Theorem
\ref{thm-mainlimit} with a table; see Table~\ref{table:r21}.\vadjust{\goodbreak}

\begin{figure}

\includegraphics{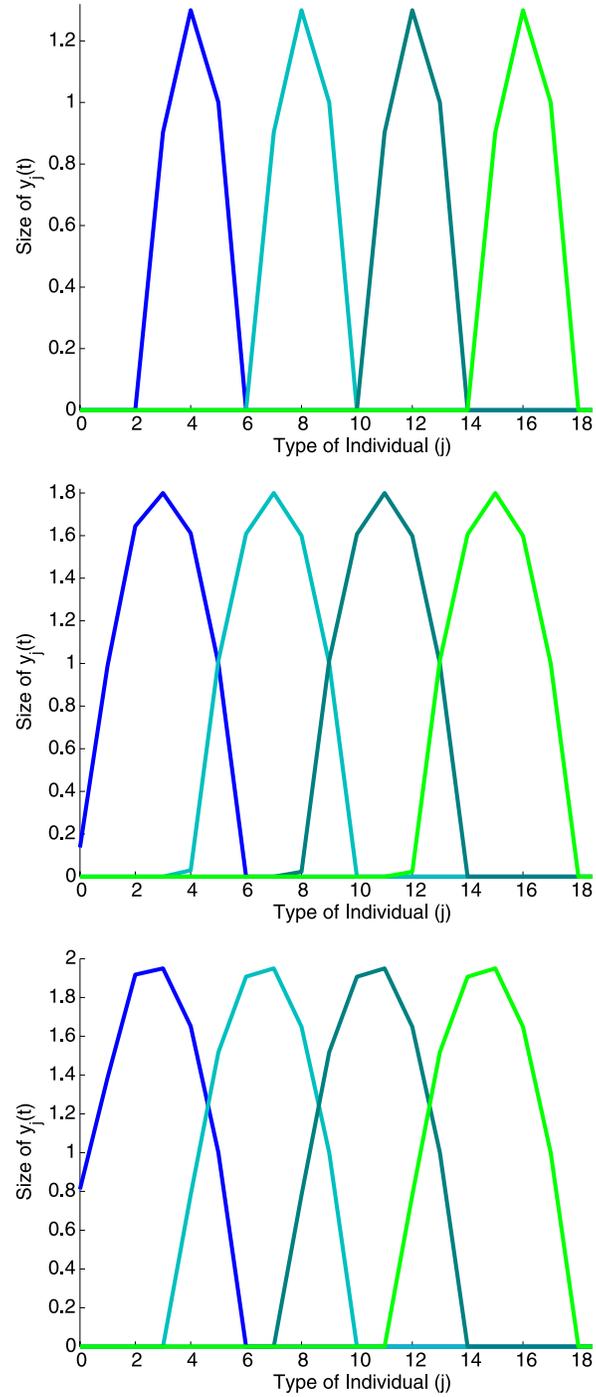}

  \caption{Distribution of types at the times $t_{4}$, $t_{8}$, $t_{12}$ and $t_{16}$
   (from left to right) with the same parameters as Figure \protect\ref{fig:r123} (top to bottom).}
  \label{fig:typedist}
\end{figure}

Since $\alpha<2$ the first two rows are the same as in regime 1, and
we again have
\[
\tau^\mu_1 \sim(2-\alpha) L/\gamma.
\]
However, we now have $\alpha> r_2$ so that
\[
\Delta_1 = \delta_{1,1} \wedge\delta_{1,2} = (\alpha-1) \wedge
(\gamma
/\gamma_2) = \gamma/\gamma_2,
\]
and hence the 2's reach level $1/\mu$ before the 1's fixate. This yields
\[
\tau^\mu_2 \sim\frac{\gamma}{\gamma_2} \cdot\frac{L}{\gamma}.
\]
Now $y_1(s_2) = 1+\gamma/\gamma_2 = r_2$, so the additional time it takes
$y_1$ to reach level $\alpha$ is $\delta_{2,1} = \alpha-r_2$. Since
$\alpha< r_3$, we have $(\alpha-r_2)\gamma_3/\gamma< 1$ and hence
$\delta_{2,1} < \delta_{2,3}$, that is, the 1's will fixate before
the 3's reach level $1/\mu$. To show that the 1's fixate before the
2's and conclude that $\Delta_2 = \delta_{2,1} = \alpha-r_2$, we
need to show that $(\alpha-r_2) < (\alpha-1) \gamma/\gamma_2, $
which holds
if and only if
%
\begin{equation}\label{keybndr2}
\alpha< \frac{2+\gamma}{1+\gamma}.
\end{equation}
However, comparing $r_3 =1 +\gamma/\gamma_2 + \gamma/\gamma_3$ with
the upper bound
in \eqref{keybndr2}, we can see that
\begin{eqnarray*}
1 +\gamma/\gamma_2 + \gamma/\gamma_3 < \frac{2+\gamma}{1+\gamma}
\hspace*{-3pt}\quad&\Longleftrightarrow&\quad\hspace*{-3pt}\frac{(3+3\gamma+\gamma^2)+(2+\gamma)}{(2+\gamma)(3+3\gamma
+\gamma^2)} < \frac{1}{1+\gamma}
\\
&\Longleftrightarrow&\quad\hspace*{-3pt}\frac{5+9\gamma+5\gamma^2+\gamma^3}{6+9\gamma
+5\gamma^2 + \gamma^3}
< 1.
\end{eqnarray*}
The last inequality is always true and therefore \eqref{keybndr2}
holds throughout regime 2 and $\Delta_2 = \alpha-r_2$, justifying the
fourth line in Table \ref{table:r21}. Finally, to check that the 2's
have not yet fixated when the 3's reach level $1/\mu$ and prove
\begin{table}[t]
\caption{$\log_{1/\mu}$ sizes in regime 2, time in units of
$L/\gamma$}
\label{table:r21}
\begin{tabular*}{277pt}{@{\extracolsep{4in minus 4in}}lcccc@{}}
\hline
\textbf{Time} & \textbf{Time increment} & \textbf{Type 1} & \textbf{Type 2} & \textbf{Type 3} \\
\hline
$0+$ & & $\alpha-1$ \\
$s_1$ & $\Delta_0 = 2-\alpha$ & 1 & 0 & \\
$s_2$ & $\Delta_1 = \gamma/\gamma_2$ & $r_2$ & 1 & 0 \\
$s_3$ & $\Delta_2 = \alpha-r_2$ & $\alpha$ & $1+ \Delta_2\gamma
_2/\gamma$
& $\Delta_2\gamma_3/\gamma$ \\
$s_4$ & $\Delta_3 = \frac{\gamma}{\gamma_2}( 1 - \Delta
_2\frac{\gamma
_3}{\gamma}) $
& & $1+ \Delta_2\gamma_2/\gamma+ \Delta_3$ & 1\\
\hline
\end{tabular*}
\end{table}
\[
\Delta_3 = \delta_{3,3} = \frac{\gamma}{\gamma_2}\biggl( 1 -
\Delta
_2\frac{\gamma_3}{\gamma}\biggr),
\]
we note that the size of $y_2(s_3 + \delta_{3,3})$ is
\begin{eqnarray*}
1+ \Delta_2 \gamma_2/\gamma+ \delta_{3,3} &=& 1 + \frac{\gamma
_2}{\gamma}
(\alpha- r_2) + \frac{\gamma}{\gamma_2}
- \frac{\gamma_3}{\gamma_2} (\alpha-r_2)
\\
&=& 1+ \gamma/\gamma_2 + \ell(\alpha-r_2)
\end{eqnarray*}
with
\[
\ell\equiv\gamma_2/\gamma- \gamma_3/\gamma_2 = \frac{(2+\gamma
)^2 - (3+3\gamma+\gamma
^2)}{2+\gamma} = \frac{1+\gamma}{2+\gamma} \in(0,1),
\]
and hence $ y_2(s_3 + \delta_{3,3}) \in(r_2,\alpha)$. This justifies
the final line of Table \ref{table:r21}, and we conclude that
\[
\frac{\tau_3^\mu}{L/\gamma} \to\Delta_2 + \Delta_3
= \alpha-r_2 + \frac{\gamma(1- \gamma_3(\alpha-r_2)/\gamma
)}{\gamma_2}.
\]

In contrast to regime 1, the relative sizes of types when the 3's reach
$1/\mu$ are not exactly the same as the relative sizes when the 2's
reach level $1/\mu$. To describe this more complicated situation,
suppose that type-$(k-2)$ individuals have size $(1/\mu)^x$ at the
time type-$(k-1)$ individuals reach level $1/\mu$. Then, if we assume:

\begin{longlist}[(2a)]
\item[(2a)] type $k-2$ reaches fixation before type $k-1$;

\item[(2b)] type $k-2$ reaches fixation before $k$'s reach $1/\mu$;

\item[(2c)] type $k$ reaches level $1/\mu$ before type $k-1$ reaches fixation,
\end{longlist}
we can repeat the arithmetic leading to Table \ref{table:r21} to yield
Table \ref{table:r22},
where here $f(x) = 1 + \gamma_2 t_k^1/\gamma+ t_k^2 = r_2 + \ell
(\alpha
-x)$ with $\ell= (1+\gamma)/(2+\gamma),$ as before.

\begin{table}[t]
\tablewidth=11,5cm
\caption{Iteration in regime 2}
\label{table:r22}
\begin{tabular*}{11,5cm}{@{\extracolsep{4in minus 4in}}lcccc@{}}
\hline
\textbf{Time increment} & \textbf{Type} $\bolds{k-3}$ & \textbf{Type} $\bolds{k-2}$ & \textbf{Type} $\bolds{k-1}$ & \textbf{Type} $\bolds{k}$ \\
\hline
& $\alpha$ & $x$ & 1 & 0\\
$t_k^1 = \alpha-x$ & $\alpha$ & $\alpha$ & $1+\gamma_2 t_k^1/\gamma
$ & $\gamma
_3 t_k^1/\gamma$ \\[2pt]
$t_k^2 = \frac{\gamma}{\gamma_2} ( 1 - \gamma_3 (\alpha-x)/\gamma
)$ & &
$\alpha$ & $f(x)$ & 1\\
\hline
\end{tabular*}
\end{table}

Since the density of 2's is $f(r_2)$ when the 3's have reached size
$1/\mu$, we see that when type $k \ge3$ reaches size $1/\mu,$ the
density of type $k-1$ is $f^{k-2}(r_2)$. This leads to the statement of
our next result.

\begin{corollary} \label{cor-threes}
Suppose $N=\mu^{-\alpha}$ for some $\alpha\in(r_2, r_3)$. Then, as
$\mu\to0,$
\begin{eqnarray*}
\frac{T_1^\mu}{L/\gamma} \to(2-\alpha), \quad\mbox{and for all $j
\geq2,$}\quad
\frac{\tau_j^\mu}{L/\gamma} \to\beta_j \qquad\mbox{in probability,}
\end{eqnarray*}
where $\beta_2=\gamma/\gamma_2,$
and if we let $f^0(x)=x,$ then for all $j \geq3$ we have
\begin{eqnarray*}
\beta_j= t_j^1 + t_j^2 =
\bigl(\alpha-f^{k-3}(r_2)\bigr) + \frac{1 - (3+3\gamma+\gamma^2) (\alpha
-f^{k-3}(r_2))}{2+\gamma}.
\end{eqnarray*}
Furthermore, the coefficients $\beta_j \to\beta_\infty$ as $j \to
\infty,$ where
\[
\beta_\infty= \alpha-r^* + \frac{ 1 - (3+3\gamma+\gamma^2)(\alpha
-r^*)}{2+\gamma}
\]
with $r^* = \lim_{j \to\infty} f^j(r_2) = (r_2 + \ell\alpha
)/(1+\ell)$.
\end{corollary}

\begin{pf}
We need to show that conditions (2a)--(2c) above are satisfied
for any $j \geq0$ and that $f^j(r_2)$ converges. The latter follows
from the fact that $f$ has slope $-\ell$ with $\ell\in(0,1)$, so, as
$j\to\infty,$
\[
f^j(r_2) \to r^* = \frac{r_2 + \ell\alpha}{1+\ell},
\]
the unique fixed point of $f$. It is easy to see that $\ell\in(0,1)$
implies that
%
\begin{equation}\label{fkbound}
r_2 \leq f^j(r_2) < \alpha
\end{equation}
for all $j \geq0,$ and (2c) immediately follows. Since $\alpha< r_3$,
we have $\gamma_3(\alpha-  r_2)/\gamma< 1$,
which, along with \eqref{fkbound}, tells us that (2b) holds for all $j
\geq0$ as well. Finally, (2a)~is equivalent to
\[
\frac{\alpha-1}{\gamma_2} > \frac{\alpha- f^{j-3}(r_2)}{\gamma}
\]
and so \eqref{fkbound} implies that to prove (2a), we need only show that
\begin{eqnarray*}
\frac{\alpha-1}{\gamma_2} > \frac{\alpha- r_2}{\gamma}.
\end{eqnarray*}
Rearranging terms, we obtain the equivalent condition $\alpha<
(2+\gamma
)/(1+\gamma)$ which holds by \eqref{keybndr2}, completing the proof.
\end{pf}

Again, the behavior of the limits $y_j(t)$ can be read from Tables \ref
{table:r21} and \ref{table:r22}.
The formulas are messy, but it is easy to compute $y_j(t)$ for a fixed
value of $\alpha$. As Figure \ref{fig:r123} shows, after a
short transient phase, the increments between the appearance of
successive types settle down into the steady-state behavior guaranteed
by Corollary~\ref{cor-threes}. Figure \ref{fig:typedist} shows the
distribution of types at various times throughout the evolution of the
system, which agree with simulations given in Figure 1 in the Appendix
of~\cite{BN}.

\subsection{Regime $j$, $j \geq3$} \label{sec-r3}

Regime 3 occurs for $\alpha\in(r_3,r_4)$ with $r_4 = r_3 + \gamma
/\gamma_4$. When $\gamma$ is small,
$\gamma/\gamma_4 \approx1/4$, so this regime is roughly $\alpha\in
(11/6,25/12)$. If $\alpha<2,$
then the initial phases are similar to regime 2, but now type 3 reaches
$1/\mu$ before the 1's fixate; see Table \ref{table:r31}.

\begin{table}[b]
\caption{$\log_{1/\mu}$ sizes in regime 3, time in units of
$L/\gamma$}
\label{table:r31}
\begin{tabular*}{225pt}{@{\extracolsep{4in minus 4in}}lcccc@{}}
\hline
\textbf{Time} & \textbf{Time increment} & \textbf{Type 1} & \textbf{Type 2} & \textbf{Type 3} \\
\hline
$0+$ & & $\alpha-1$ \\
$s_1$ & $\Delta_0 = 2-\alpha$ & 1 & 0 & \\
$s_2$ & $\Delta_1 = \gamma/\gamma_2$ & $r_2$ & 1 & 0 \\
$s_3$ & $\Delta_2 = \gamma/\gamma_3$ & $r_3$ & $1+ \gamma_2/\gamma
_3$ & 1
\\
\hline
\end{tabular*}
\end{table}

\begin{figure}[b] 

\includegraphics{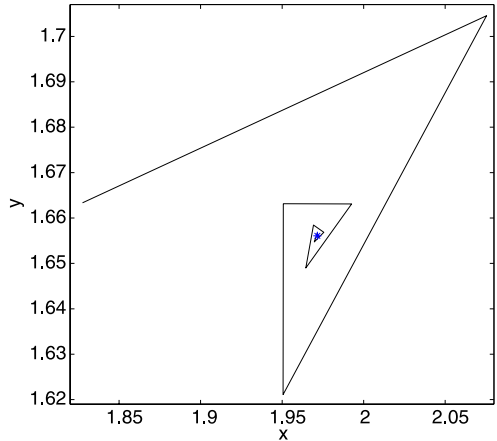}

  \caption{Successive iterates of the two-dimensional map $(x,y) \mapsto f(x,y) = (f_1(x,y),f_2(x,y))$
  given in Table \protect\ref{table:r3}, started with initial conditions $(x_0,y_0) = (1+\gamma/\gamma_2+ \gamma/\gamma_3, 1 + \gamma_2/\gamma_3)$
  and parameters $\alpha = 1.95$, $\gamma = 0.01$. The star denotes the fixed point of $f$.}
  \label{fig:r3its}
\end{figure}

If we now assume that:
\begin{longlist}
\item[(3a)] type $k-3$ reaches fixation before types $k-2$ and $k-1$;

\item[(3b)] type $k-3$ reaches fixation before type $k$'s reach $1/\mu$;

\item[(3c)] type $k$ reaches level $1/\mu$ before types $k-2$ and $k-1$
reaches fixation,
\end{longlist}
then the recursion in Table \ref{table:r22} becomes a pair of
equations (see Table \ref{table:r3}). To imitate the proof in regime 2
we would have to show that (3a)--(3c) hold for $x=r_3$ and
$y=1+\gamma_2/\gamma$, and for all of the iterates $f^k(x,y),$ where $f
\equiv(f_1,f_2)$. Figure~\ref{fig:r3its} shows that this is true when
$\alpha=1.95$ and $\gamma=0.01$; however, verifying this
algebraically is difficult because $f(x,y)$ may fail to satisfy the
conditions when $(x,y)$ does.

\begin{table}[t]
\caption{Iteration in regime 3}
\label{table:r3}
\begin{tabular}{lccccc}
\hline
\textbf{Time increment} & \textbf{Type} $\bolds{k-4}$ & \textbf{Type} $\bolds{k-3}$ & \textbf{Type} $\bolds{k-2}$ & \textbf{Type} $\bolds{k-1}$ &
\textbf{Type} $\bolds{k}$\\
\hline
& $\alpha$ & $x$ & $y$ & 1 & \\
$t_k^1 = \alpha-x$ & $\alpha$ & $\alpha$ & $y+t_k^1\gamma_2/\gamma
$ & $1 +
t_k^1\gamma_3/\gamma$ & $t^1_k\gamma_4/\gamma$ \\[2pt]
$t_k^2 = \frac{\gamma}{\gamma_3} ( 1 - t^k_1\gamma_4/\gamma) $ & &
$\alpha$ &
$f_1(x,y)$ & $f_2(x,y)$ & 1\\
\hline
\end{tabular}
\end{table}

\begin{figure}

\includegraphics{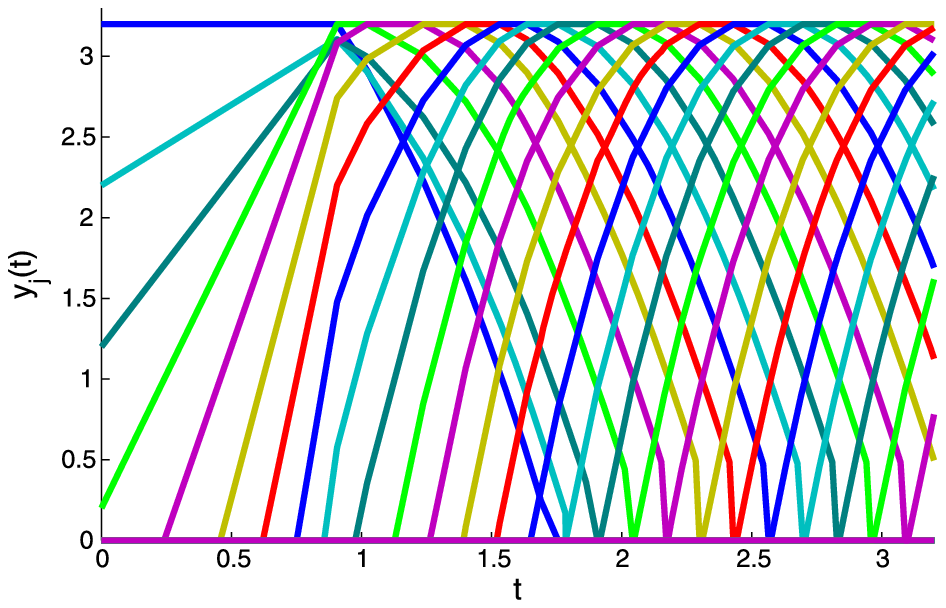}

  \caption{Plot of the limiting dynamical system in Theorem \protect\ref{thm-mainlimit}
  (fixed population size). Parameters: $\gamma = 0.1$, $\alpha=3.2$.}
  \label{fig:rinf}
\end{figure}

In general, we conjecture that if we define
\[
r_j = \sum_{i=1}^j (\gamma/\gamma_i),
\]
then we are in regime $j$ if $\alpha\in[r_j, r_{j+1})$ and we have
$\beta_i \to\beta_\infty$ as $i \to\infty$. In particular, this
would imply that $t^* = \infty$ as long as
\[
\alpha< r_\infty= \sum_{i=1}^\infty(\gamma/\gamma_j).
\]
However, as Figure \ref{fig:rinf} shows, the converse is not true.
There, we have $\alpha=3.2, \gamma=0.1$ so that $\alpha> r_\infty
\approx3.1$, but it appears that we still approach a constant
increment between waves.

\subsection{Blow-up in finite time} \label{sec-bu}

In this section we prove Lemma \ref{blowup}, which shows that for any
$\gamma>0$, we can choose $\alpha$ large enough to make $t^*$ finite. To
this end, let $\gamma> 0,$ and for all $j \geq1$ define
\[
S_j = S_j (\gamma) = \sum_{i=0}^\infty\frac{\gamma_j}{\gamma_{j+i}}
\]
and let $S = S(\gamma) = \sup\{S_j(\gamma)\dvtx j \geq1\}$. Note that since
\[
S_j = \sum_{i=0}^\infty\frac{1 - (1+\gamma)^{-j}}{(1+\gamma)^i -
(1+\gamma
)^{-j}} \le\sum_{i=0}^\infty\frac{1}{(1+\gamma)^i} \cdot\frac
{1}{1-(1+\gamma)^{-j-i}},
\]
we have $S < \infty$. With this notation to hand, we can prove
the lemma.

\begin{pf*}{Proof of Lemma \ref{blowup}}
Fix $\gamma>0$, define $S = S(\gamma)$ as above and choose $\alpha$ large
enough so that $\alpha> 1+2S$ and
\[
\gamma_{a/2}/\gamma_{a} < 1/S,
\]
where $a \equiv\lfloor\alpha\rfloor$. We will show that $\Delta_1
= (1-(\alpha-a))\gamma/\gamma_a$ and $\Delta_n \leq\gamma/\gamma
_{a+n-1}$ for $n
\geq2$ so that
\[
t^* = \sum_{n=1}^\infty\Delta_n \leq\gamma\sum_{n=0}^\infty
[(1+\gamma
)^{a+n} -1]^{-1} < \infty.
\]
To prove that $\Delta_n$ has the desired bound, we will show that for
all $n \geq0$, $y_{a+n}(t)$ hits level 1 before $y_j(t)$ hits level
$\alpha$ for any $1 \leq j \leq a+n$. Since this implies that the
growth rate of type $(a+n)$'s is $\gamma_{a+n}/\gamma$, we have, in the
notation of Theorem~\ref{thm-mainlimit}, $k_n^* = a+n-1$ for all $n
\geq1$ so that $\Delta_1 = \delta_{1,k_1^*} = (1-(\alpha-a))\gamma
/\gamma
_a$ and $ \Delta_n = \delta_{n,k_n^*} = \gamma/\gamma_{a+n-1}$ for
all $n
\geq2$.

We first note that since $y_j(0) = (\alpha-j)^+$ for all $j \geq0$,
we know that no type $1 \leq j \leq a-1$ can reach level $\alpha$
before type-$(a-1)$ individuals reach level 1. Now, let $n \geq0$ and
suppose that no individual of type $1 \leq j \leq a+(n-1)$ reaches
level $\alpha$ before type-$(a+(n-1))$ individuals reach level 1.
Then, for $a+1 \leq j \leq a+n$, type-$(j-1)$ individuals reach level 1
at time
\[
t_{j-1} = \bigl(1-(\alpha-a)\bigr)\gamma/\gamma_a + \sum_{i=1}^{j-a-1} \frac
{\gamma}{\gamma_{a+i}},
\]
at which time type-$j$ individuals are born and start to grow at rate
$\gamma_{j}/\gamma$. If there is no change in the dominant type, then type
$j$'s will reach level $\alpha$ at time $j \gamma/\gamma_j$ if $j
\leq a$
and time $t_{j-1} + \alpha\gamma/\gamma_j$ if $j > a,$ so if we
define $t_j
= 0$ for all $j \leq a$, then the proof will be complete if we can show that
%
\begin{equation}\label{bndj}
\min(j,\alpha) \gamma/\gamma_j > (t_{a+n-1}-t_{j-1}) + \gamma
/\gamma_{a+n}
\end{equation}
for all $1 \leq j \leq a+n$. Suppose first that $1 \leq j < a/2$. Our
choice of $\alpha$ then implies that
\[
j \geq1 > \frac{\gamma_{a/2}}{\gamma_a}   S> \frac{\gamma
_j}{\gamma_a} \sum
_{i=0}^{n} \frac{\gamma_a}{\gamma_{a+i}}
\]
so that
\[
j(\gamma/\gamma_j) > \sum_{i=0}^{n} \frac{\gamma}{\gamma_{a+i}} >
t_{a+n-1} +
\frac{\gamma}{\gamma_{a+n}}.
\]
If $a/2 \leq j \leq a$, we have
\[
j \geq a/2 > S > \sum_{i=0}^n \frac{\gamma_a}{\gamma_{a+i}} > \sum_{i=0}^n
\frac{\gamma_j}{\gamma_{a+i}}
\]
so that, again,
\[
j(\gamma/\gamma_j) > t_{a+n-1} + \frac{\gamma}{\gamma_{a+n}}.
\]
Finally, if $ a+1 \leq j \leq a+n,$ we have
\[
\alpha> S > \sum_{i=0}^\infty\frac{\gamma_j}{\gamma_{j+i}} >
(\gamma_j/\gamma
)(t_{a+n-1}-t_{j-1} + \gamma/\gamma_{a+n}),
\]
which completes the proof of \eqref{bndj}.
\end{pf*}

\section{Ideas behind the proof} \label{pf-ideas}
For the remainder of the paper, we suppose that $\rho\geq0$ and that
$(\alpha,\gamma) \in G_\rho$. $C$ will always denote a constant that
does not depend on $\mu$ and whose value may change from line to line.
We begin with a simple, but useful, lemma which explains why the
limiting result for the birth times follows from the limiting result
for the sizes.

\begin{lemma} \label{nonewguys}
Let $\varepsilon, b >0$. Then,
\[
P\bigl(T_{j+1}^\mu\leq Lt,   X_j^\mu(Ls) \leq b   (1/\mu
)^{1-\varepsilon
} \mbox{ for all } s \leq t\bigr) \to0
\]
as $\mu\to0$ for any $j \geq0$.
\end{lemma}

\begin{pf}
Let $M_j(t)$ denote the number of mutations from $j$'s to $(j+1)$'s by
time $t$ and let
\[
A(t) = \{X_j^\mu(Ls) \leq b (1/\mu)^{1-\varepsilon} \mbox{ for all
} s \leq
t\}.
\]
Since mutations to $(j+1)$'s occur at rate $\mu X_j^\mu(t),$ we have
\[
E   [M_j(Lt); A(t)] \leq b Lt\mu^\varepsilon\to0
\]
as $\mu\to0$ and therefore Chebyshev's inequality implies that
\begin{eqnarray*}
P\bigl(M_j(Lt) \geq1, A(t)\bigr) \leq E  [ M_j(Lt); A(t)] \to0
\end{eqnarray*}
as $\mu\to0,$ yielding the result.
\end{pf}

Assuming we have the uniform convergence of $Y_i^\mu(t) \to y_i(t)$
for all $i \geq0$, Lemma \ref{nonewguys} implies that $P(T_{j+1}^\mu
\leq(t_j-\varepsilon)L/\gamma) \to0$, but since $Y_{j+1}^\mu(t)
\to
y_{j+1}(t)$ also implies that
\[
P\bigl(T_{j+1}^\mu> (t_j+\varepsilon)L/\gamma\bigr) \leq P\bigl(X_{j+1}^\mu
\bigl((t_j+\varepsilon)L/\gamma\bigr) =
0\bigr) \to0,
\]
the desired convergence of $T_{j+1}^\mu$ follows.

Our next result gives an approximation for the population size $N^\mu
(t)$ that yields the desired uniform convergence of $F^\mu(t)$ and
also proves useful in other situations.

\begin{lemma} \label{poplim}
Let $\zeta,a >0$. Then, as $\mu\to0$,
\begin{eqnarray*}
P\biggl(\sup_{0\leq t \leq a L} \biggl|\frac{N^\mu(t)}{N^\mu(0)
e^{\rho t}} - 1\biggr| > \zeta\biggr) \to0.
\end{eqnarray*}
\end{lemma}

\begin{pf}
Let $N_j(t)$, $0 \leq t \leq a L,$ be a family of i.i.d.~pure birth
(Yule) processes in which individuals give birth at rate $\rho$ and
the initial population is $N_j(0)=1$. We then have
\begin{eqnarray*}
N^\mu(t) \stackrel{d}{=} \sum_{j=1}^{N^\mu(0)} N_j(t).
\end{eqnarray*}
It follows, for example, from \cite{AN}, page 109, equation (5), that
the moments $m_j^i(t) = E (N_j^i(t))$, $i=1,2,$ satisfy
\begin{eqnarray*}
m_j^1(t) &=& e^{\rho t},
\\
m_j^2(t) &=& 2 e^{2\rho t}(1- e^{-\rho t}) \leq C e^{2 \rho t},
\end{eqnarray*}
and so $M_j(t) = e^{-\rho t}N_j(t)-1$, $t \geq0$, is a mean-zero
martingale (\cite{AN}, page 111) with
\[
\operatorname{var}(M_j(t)) = \frac{m_j^2(t)}{e^{2\rho t}} - 1 \leq C.
\]
Since the $M_j$ are independent,
\begin{eqnarray*}
M(t) = \sum_{j=1}^{N^\mu(0)} M_j(t)
\end{eqnarray*}
(which is itself a mean-zero martingale) has
\[
E(M_j^2(t)) = \operatorname{var}(M_j(t)) \leq C N^\mu(0).
\]
Applying Chebyshev's inequality and the $L^2$-maximal inequality yields\vspace*{2pt}
\begin{eqnarray*}
&&P\biggl(\sup_{0\leq t \leq a L}
\biggl|\frac{\sum_{j=1}^{N^\mu(0)} N_j(t)}
{N^\mu(0) e^{\rho t}} - 1\biggr| > \zeta\biggr)
\\[2pt]
&&\qquad\leq P\Bigl(\sup_{0 \leq t \leq a L} |M(t)| > \zeta N^\mu(0)\Bigr)
\\[2pt]
&&\qquad\leq\frac{4}{\zeta^2 N^\mu(0)^2 } E (M^2(aL))
\\[2pt]
&&\qquad\leq\frac{C}{\zeta^2 N^\mu(0)} \to0 \qquad\mbox{as $\mu\to0,$}\vspace*{2pt}
\end{eqnarray*}
which gives the desired result.
\end{pf}

There are four steps involved in proving the desired convergence of
$Y_j(t)$, \mbox{$ j \geq0$}, in Theorem \ref{induct-grow}. The first step,
taken in Section \ref{sec-startup}, is to prove a result about the
initial behavior of the process.

\begin{proposition}\label{startup}
Let $k=\lfloor\alpha\rfloor$ be the largest integer $\leq\alpha$
and define
\[
\delta_{0,j} =
\cases{
j \gamma/(\lambda_j-\rho), &\quad $j < k$,
\vspace*{2pt}\cr
\bigl(1-(\alpha-k)\bigr)\gamma/\lambda_k, &\quad $j=k$.
}\vadjust{\goodbreak}
\]

\noindent
Then, for any $0 < t_1 < t_2 < \Delta_0 \equiv\min\{\delta_{0,j}\dvtx j
\leq k\}$, $Y_j^\mu(t) \to y_j(t)$ in probability uniformly on
$[t_1,t_2]$ with
\begin{eqnarray*}
y_j(t) =
\cases{ (\alpha-j) + t \lambda_j/\gamma, &\quad $j \leq k$, \cr
0, &\quad $j > k$.
}
\end{eqnarray*}
\end{proposition}

Proposition \ref{startup} yields the correct initial conditions (a).
The proof of the inductive step (b) is given in Section \ref
{sec-inductpf} and has three main parts that together roughly describe
how the limit changes during one iteration of (b), that is, on the
interval $[s_n,s_{n+1}]$. Since we wish to apply the results below to
$Y_j^\mu(t)$ at positive times, we consider a version of our Moran
model in which we allow for general initial conditions $X^\mu(0)$
satisfying the following.

\begin{Assumptions*}  As $\mu\to0$, $F^\mu(0) \to\alpha>0$ and
$Y^\mu_j(0) \to y_j^0$ in probability for all $j \geq0$. Furthermore,
we suppose that the $y_j^0,   j \geq0,$ satisfy the conditions:
\begin{longlist}[(iii)]
\item[(i)] there is a unique value of $m$ with $y_m^0= \alpha$;

\item[(ii)] there is a $k > 0$ such that $y_j^0=0$ for all $j>k$, $y_j^0>0$
for $m<j\le k$ and $y_k^0<1$;

\item[(iii)] $y_{j+1}^0 > y_j^0 - 1$ for $0 \le j \le k$.
\end{longlist}
\end{Assumptions*}

Define
\[
\delta^\varepsilon_j \equiv
\cases{ (\alpha-y_j^0-\varepsilon)\gamma/(\lambda_{j-m}-\rho
), &\quad $m < j <
k$, \cr
(1- y_k^0-\varepsilon)\gamma/\lambda_{k-m,} &\quad $j=k$
}
\]
and let $\Delta_\varepsilon\equiv\min\{\delta^\varepsilon_j\dvtx m <
j \le k\}$ for
$\varepsilon\geq0$. For $j \geq0$ and $t \leq\Delta_0$, define
\begin{eqnarray*}
y_j(t) =
\cases{ (y_j^0 + t\lambda_{j-m}/\gamma)^+, &\quad $j \leq k$,
\cr
0, &\quad $j > k$.
}
\end{eqnarray*}
To connect the next three results below back to (b), we will use
Proposition \ref{interior} to describe the limit on the intervals
$[s_n+\varepsilon,s_{n+1}-\varepsilon']$ for small $\varepsilon,
\varepsilon' >0$ and use
Propositions \ref{case1} and \ref{case2} to describe the limit on
$[s_{n+1}-\varepsilon',s_{n+1}+\varepsilon],$ depending on which of
the following two
possible outcomes occurs: (i) $\Delta_0 = \delta^0_k$ and a new type
is born, or (ii) $\Delta_0 = \delta^0_{n} $ for some $n \in(m,k)$
and there is a change in the dominant type.

\begin{proposition} \label{interior}
Let $\varepsilon>0$ and suppose that the above assumptions hold. Then,
$Y^\mu
_j(t) \to y_j(t)$ in probability uniformly on $[0,\Delta_\varepsilon
]$ for
all $j \geq0$.
\end{proposition}

\begin{proposition} \label{case1}
Suppose that the above assumptions hold and that $\Delta_0 = \delta
_k^0$. For $t \leq\varepsilon$, let
\begin{eqnarray*}
y_j(\Delta_0 + t) =
\cases{ \bigl(y_j(\Delta_0) + t\lambda_{j-m}/\gamma\bigr)^+, &\quad $j \le
k+1$, \cr
0, &\quad $j > k+1$.
}
\end{eqnarray*}
There then exists $\varepsilon_1 = \varepsilon_1(y^0) > 0$ such that
for all $j \neq
k+1$, $Y^\mu_j(t) \to y_j(t)$ in probability uniformly on $[\Delta
_{\varepsilon},\Delta_0+\varepsilon]$ and
%
\begin{eqnarray} \label{kplus1}
&& P\Bigl( \sup_{ \Delta_{\varepsilon/2} \leq t \leq\Delta
_0+\varepsilon}
Y_{k+1}(t) - (t-\Delta_{\varepsilon/2}) \lambda_{k+1-m}/\gamma>
\varepsilon/2\Bigr)
\to0, \nonumber
\\
&&\hspace*{7,7pt}P\Bigl( \sup_{ \Delta_0+(\gamma/\lambda_{k-m}) \varepsilon/2
\leq t \leq
\Delta_0+\varepsilon} Y_{k+1}(t)
\\
&&\hspace*{7,7pt}\qquad{}- \bigl(t-\Delta_0-(\gamma/\lambda
_{k-m}) \varepsilon
/2\bigr)\lambda_{k+1-m}/\gamma< -\varepsilon/2\Bigr) \to0\nonumber
\end{eqnarray}
as $\mu\to0,$ provided $\varepsilon< \varepsilon_1$.
\end{proposition}

\begin{proposition}\label{case2}
Suppose that the above assumptions hold and that $\Delta_0 = \delta
_n^0$ for some $ n \in(m,k)$. For $t \le\varepsilon$, let
\begin{eqnarray*}
y_j(\Delta_0 + t) =
\cases{ \bigl(y_j(\Delta_0) + t\lambda_{j-n}/\gamma\bigr)^+, &\quad $j \le k$,
\cr
0, &\quad $j > k$.
}
\end{eqnarray*}
There then exists $\varepsilon_2 = \varepsilon_2(y^0) > 0$ such that
$Y^\mu_j(t) \to
y_j(t)$ in probability uniformly on $[\Delta_{\varepsilon},\Delta
_0+\varepsilon],$
provided $\varepsilon< \varepsilon_2$.
\end{proposition}

 Note that $n \neq m+1$ is possible (see Figure \ref{fig:rinf}).

\begin{pf*}{Proof of Theorem \ref{induct-grow} from Propositions
\ref{startup}--\ref{case2}} Suppose that $X_0^\mu(0) = N^\mu(0)$ and
$X_j^\mu(0) = 0$ for all $j \geq1$, and let $y_j(t)$ denote the
dynamical systems described by (a) and (b). Let $K$ be a compact subset
of $(0,t^*)$, $\zeta>0$ and take $a \in(0,\Delta_0)$, $n(K) \geq1$
so that $[a,s_{n(K)}]\supset K$, where $s_n$ is as defined in (b).
Choose $\varepsilon>0$ small enough so that $\varepsilon< \varepsilon
_1(y(s_n)), \varepsilon
_2(y(s_n))$ for all $n \leq n(K),$ where $\varepsilon_1,\varepsilon
_2$ are as in
Propositions \ref{case1} and \ref{case2}, respectively. Without loss
of generality, suppose that $\varepsilon< \zeta/(c+1),$ where $c =
c(\gamma,\rho
) >1$ is defined below. We also set $s_{n,\varepsilon} = s_n -
\varepsilon\gamma/\lambda
_{j^*_{n-1}-m_{n-1}}$, where $j^*_n$ satisfies $\Delta_{n} = \delta
_{n,j^*_{n}}$.

By Proposition \ref{startup}, we obtain $Y_j(t) \to y_j(t)$ in
probability uniformly on $[a,s_{1,\varepsilon}]$. Suppose now that we have
uniform convergence on $[a,s_{n,\varepsilon}]$ for some \mbox{$n \leq
n(K)-1$}. We
then have two cases to consider. If $j^*_n = j$ for some $j \in
(m_n,k_n)$, then applying Proposition \ref{case2} up to time
$s_n+\varepsilon
$ and then Proposition \ref{interior} with $y_j^0 =
y_j(s_n+\varepsilon)$ for
all $j$ up to time $s_{n+1,\varepsilon}$, we obtain the result. If
$j^*_n =
k_n$, then Proposition \ref{case2} clearly allows us to extend uniform
convergence for $Y_j(t)$, $j \neq k_n+1$ up to time $s_n+\varepsilon$.
To do
this for $j=k_n+1$, we first apply Proposition~\ref{interior} to get
convergence up to time $s_{n,\varepsilon/2}$. Write
\begin{eqnarray*}
Y_j(t) - y_j(t) &=&\bigl(Y_j(t) - (t-s_{n,\varepsilon/2}) \lambda
_{k_n+1-m_n}/\gamma\bigr)
 \\
 &&{}+
\bigl( (t-s_{n,\varepsilon/2}) \lambda_{k_n+1-m_n}/\gamma-y_j(t)\bigr).
\end{eqnarray*}
Recalling that $y_j(t) = 0$ if $t \leq s_n$ and $y_j(t) =(t-s_n)
\lambda_{k_n+1-m_n}/\gamma$ if $s_n \leq t \leq s_n+ \varepsilon$,
we can see that
\begin{eqnarray*}
(t-s_{n,\varepsilon/2}) \lambda_{k_n+1-m_n}/\gamma-y_j(t) \in
[0,(\lambda
_{k_n+1-m_n}/\lambda_{k_n-m_n})\varepsilon/2] \subset[0,c
\varepsilon/2]
\end{eqnarray*}
for all $s_{n,\varepsilon/2} \leq t \leq s_n+\varepsilon$, the last inclusion
following from the fact that
\begin{eqnarray*}
\lambda_{k+1}/\lambda_k &=& \bigl((1+\rho)(1+\gamma)^{k+1}-1\bigr)/\bigl((1+\rho
)(1+\gamma
)^k-1\bigr)
\\
&\leq&\bigl((1+\rho)(1+\gamma)^2-1\bigr)/\bigl((1+\rho)(1+\gamma)-1\bigr)
\\
&\equiv& c
\end{eqnarray*}
for all $k\geq1$. Since Proposition \ref{case1} implies that $Y_j(t)
- (t-s_{n,\varepsilon/2})\lambda_{k_n+1-m_n}/\gamma< \varepsilon/2$
for all $s_{n,\varepsilon
/2} \leq t \leq s_n+\varepsilon$ with high probability and $c > 1$, we obtain
\begin{eqnarray*}
P\Bigl(\sup_{s_{n,\varepsilon/2} \leq t \leq s_n+\varepsilon}
Y_{k+1}^\mu(t) -
y_{k+1}(t) > (c+1) \varepsilon/2 \Bigr) \to0
\end{eqnarray*}
as $\mu\to0$. To prove the lower bound, we note that $Y_{k+1}(t) -
y_{k+1}(t) \geq0$ for $t \leq s_n$, $y_{k+1}(t) \leq c \varepsilon/2$
for all
$t \leq s_n + (\gamma/\lambda_{k_n-m_n}) \varepsilon/2$ and, by a similar
argument to the one above, using the second equation in \eqref{kplus1}
instead of the first, $Y_{k+1}(t) - y_{k+1}(t) < -(c+1)\varepsilon/2$
for all
$s_n + (\varepsilon/2)(\gamma/\lambda_{k_n-m_n}) \leq t \leq s_n +
\varepsilon$ with
high probability. Therefore,
\begin{eqnarray*}
P\Bigl(\sup_{s_{n,\varepsilon/2} \leq t \leq s_n+\varepsilon}
Y_{k+1}^\mu(t) -
y_{k+1}(t) < - (c+1) \varepsilon/2 \Bigr) \to0.
\end{eqnarray*}
Since $\varepsilon< \zeta/(c+1) $, we conclude that
\begin{eqnarray*}
P\Bigl(\sup_{s_{n,\varepsilon/2} \leq t \leq s_n+\varepsilon}
|Y_{k+1}^\mu(t) -
y_{k+1}(t)| > \zeta\Bigr) \to0
\end{eqnarray*}
as $\mu\to0,$ so we have convergence up to time $s_n+\varepsilon$. Finally,
to complete the proof of the inductive step, apply Proposition \ref
{interior} with $y_j^0 = y_j(s_n+\varepsilon)$ to extend the
convergence up to
time $s_{n+1,\varepsilon}$.
\end{pf*}

\section{Initial behavior} \label{sec-startup}

In this section we prove Proposition \ref{startup} concerning the
initial behavior of the limit, but before we can begin, we need to take
care of some preparatory details. We set $\mathbb{N}_{0} = \{
0,1,\ldots\}$ and for $x = (x_0,x_1,\ldots ) \in\mathbb{R}^{\mathbb
{N}_0}$, we write $x^{j,k} = x + e_j - e_k$, where the $e_j \in\mathbb{R}
^{\mathbb{N}_0}$, $j \geq0$, are the standard basis vectors. It is
useful to note that we can define $\{(N^\mu(t),X^\mu(t))\}_{t \geq
0}$ as the Markov process with state space
\begin{eqnarray*}
\mathcal{S} \equiv\biggl\{(N,x) \in\mathbb{N}_0 \times\mathbb
{N}_0^{\mathbb{N}_0}\dvtx \sum_{j\geq0} x_j = N\biggr\}
\end{eqnarray*}
and initial population $(N^\mu(0),X^\mu(0))= (N^\mu(0),(N^\mu
(0),0,0,\ldots))$ with $N^\mu(0)$ distributed according to $\nu_0$ in
which $(N,x) \mapsto(N,y)$ at rate $p_{j,k}(x) + \mu\delta_{j-1,k}
x_{j-1}\hspace*{-1pt}$ if $y = x^{j,k}$ for some $j,k \geq0$, $(N,x) \mapsto
(N+1,y)$ at rate $\rho N (1+\gamma)^j x_j/w $ if $y = x + e_j$ and $(N,x)
\mapsto(M,y)$ at rate 0 otherwise, where $\delta_{j,k}$ here denotes
the Kronecker delta symbol and
\begin{eqnarray*}
p_{j,k}(x) = \frac{(1+\gamma)^j x_j x_k}{w}, \qquad w = \sum_{i \geq0}
(1+\gamma)^i x_i.
\end{eqnarray*}
We let
\begin{eqnarray*}
b^0_j(x) = \rho N (1+\gamma)^j x_j/w + \sum_{k \neq j} p_{j,k}(x),
\qquad
d^0_j(x) = \sum_{k \neq j} p_{j,k}(x)
\end{eqnarray*}
denote the birth and death rates, respectively, of type $j$'s, ignoring
mutations, and drop the 0's when the mutation rates are included.
$\mathcal{F}_t = \sigma\{X^\mu(s)\dvtx s \leq t\},$ and unless otherwise
explicitly stated, when we say a process is a martingale,
submartingale, etc., it will be with respect to the canonical
filtration $\mathcal{F}_t$. We will also use the notation
\[
\mathcal{S}^N = \biggl\{x \in\mathbb{N}_0^{\mathbb{N}_0}\dvtx \sum_{j\geq0}
x_j = N\biggr\}
\]
to denote a particular cross section of our state space $\mathcal{S}$.

For convenience, we will assume for the remainder of this section that
$N^\mu(0) = \mu^{-\alpha}$. Our first lemma, which is similar in
spirit to Lemma \ref{nonewguys}, takes care of the limits for $j \geq
k+1$. Recall that $T^\mu_{k+1} = \min\{ t \dvtx X^\mu_{k+1}(t)>0 \}$.

\begin{lemma} \label{nobettertype}
If $k=\lfloor\alpha\rfloor,$ then
$P(T_{k+1}^\mu< Lt /\gamma) \to0$ as $\mu\to0$ for any $t < \delta_{0,k}$.
\end{lemma}

\begin{pf}
Since type $j$'s are born at rate $b_j(x)$ and die at rate $d_j(x)$, we have
\begin{eqnarray*}
\frac{d}{dt}E X_j^\mu(t) = E\bigl(b_j(X^\mu(t))- d_j(X^\mu(t))\bigr).
\end{eqnarray*}
Using $\sum X_i^\mu(t) = N^\mu(t)$ and $(1+\gamma)^i \ge1$ for $i
\geq0$, we have
%
\begin{eqnarray} \label{lambda}
b_j(X^\mu(t))-d_j(X^\mu(t)) &=& \frac{\sum_{i \geq0} [(1+\rho
)(1+\gamma)^j-(1+\gamma)^i]X_i^\mu(t) X_j^\mu(t) }{\sum_{i \geq0}
(1+\gamma
)^i X_i^\mu(t)} \nonumber
\\
&&{}+ \mu\bigl(X_{j-1}^\mu(t)-X_j^\mu(t)\bigr)
\\
&\leq&\lambda_j X_j^\mu(t) + \mu X_{j-1}^\mu(t)\nonumber
\end{eqnarray}
for any $t \geq0$. Thus, for $j \geq1$, we obtain
\begin{eqnarray*}
\frac{d}{dt}E X_j^\mu(t) \leq\lambda_j E X_j^\mu(t) + \mu E
X_{j-1}^\mu(t)
\end{eqnarray*}
so that integrating both sides yields
\[
E X_j^\mu(t) \leq\mu\int_0^t E X_{j-1}(s) e^{\lambda_j (t-s)}\, ds
\qquad\mbox{for $j \geq1$.}
\]
We claim that induction now implies
%
\begin{equation}\label{Xjmeanbnd}
E X_j^\mu(t) \leq C_j (1/\mu)^{\alpha-j} e^{\lambda_j t}.
\end{equation}
To prove this, we note that $E X_0^\mu(t) \leq E N^\mu(t) = (1/\mu
)^\alpha e^{\rho t}$ [recall that $N^\mu(t)$ is just a Yule process],
so the result for $j=0$ holds with $C_0=1$. Using the induction
hypothesis and integrating, we have
\begin{eqnarray*}
E X_j^\mu(t) & \leq&\mu\int_0^t C_{j-1} (1/\mu)^{\alpha
-j+1}e^{\lambda_{j-1}s} e^{\lambda_j(t-s)}   \,ds
\\
& \leq& C_{j-1} (1/\mu)^{\alpha-j} e^{\lambda_j t} \int_0^t
e^{-(\lambda_j-\lambda_{j-1})s}  \, ds,
\end{eqnarray*}
which proves the claim with $C_j = C_{j-1}/(\lambda_j-\lambda_{j-1})$.

From \eqref{Xjmeanbnd}, it follows that
\begin{eqnarray*}
\int_0^t E X_j^\mu(s)  \,ds \leq C (1/\mu)^{\alpha-j} e^{\lambda_j t}.
\end{eqnarray*}
In particular, taking $t < \delta_{0,k} = \gamma(1-(\alpha
-k))/\lambda
_k$, we have
%
\begin{equation}\label{Xkmeanbnd}
\int_0^{Lt/\gamma} E X_k^\mu(s) \, ds \leq C (1/\mu)^{1-(\delta
_{0,k}-t)\lambda_k/\gamma}.
\end{equation}
The rest of the proof is the same as the proof of Lemma \ref{nonewguys}.
\end{pf}


To obtain the appropriate limits for $j \leq k$ and complete the proof
of Proposition~\ref{startup}, we will couple $X_j^\mu(t)$, $j \leq
k,$ with upper- and lower-bounding branching processes $Z_{j,u}^\mu
(t)$ and $Z_{j,\ell}^\mu(t)$, respectively, so that $Z_{j,\ell}^\mu
(t) \leq X_j^\mu(t) \leq Z_{j,u}^\mu(t)$ up until some stopping time
$\sigma,$ which will be greater than $Lt/\gamma$ with high probability
for any $t < \Delta_0,$ and will then show that we have
\[
(1/L) \log^+ Z_{j,a}(Lt/\gamma) \to y_j(t)
\]
in probability uniformly on $[t_1,t_2]$ for any $0 < t_1 < t_2 < \Delta
_0$ (see Lemma \ref{propcompbp}). The coupling is made possible by
applying the following result to bound the birth and death rates of
type $j$'s on the interval $[0,\Delta_0]$.

\begin{lemma} \label{bdbnds}
Suppose that $x \in\mathcal{S}^N$ and that there exist $m,M \in
\mathbb{N}_0$, $\eta>0$ such that \textup{(i)} $\sum_{j \neq m} x_j \le\mu
^\eta N$ and \textup{(ii)} $x_j= 0$ for all $j > M$. For all $j \neq m$, we then
have the inequalities
\begin{eqnarray*}
\frac{(1+\rho-\mu^\eta)(1+\gamma)^{j-m}x_j}{1+ g_\mu} &\leq& b^0_j(x)
\leq\frac{(1+\rho)(1+\gamma)^{j-m}x_j}{1-h_\mu},
\\
\frac{(1-M\mu^\eta)x_j}{1+ g_\mu} &\leq& d^0_j(x) \leq x_j,
\end{eqnarray*}
where $g_\mu= \gamma_{M-m} (M-m) \mu^\eta$ and $h_\mu= -\gamma
_{-m} m \mu^\eta.$
\end{lemma}

\begin{pf} From the definition,
\begin{eqnarray*}
b_j^0(x) &=& x_j(1+\gamma)^j \frac{ (1+\rho)N - x_j }{ \sum_i
(1+\gamma)^i x_i }
\\
&=& x_j(1+\gamma)^{j-m} \frac{ (1+\rho) - x_j/N }{ 1 + \sum_i
[(1+\gamma)^{i-m} - 1 ] x_i/N }.
\end{eqnarray*}
To get the lower bound, drop the terms in the denominator with $i\le
m$, which are $\le0$,
and use the fact that $j \to\gamma_j$ is increasing.
For the upper bound, drop the terms with $i\ge m$. The death rates are
given by
\[
d_j^0(x) = x_j \frac{\sum_{i\neq j} (1+\gamma)^i x_i }{\sum_{i}
(1+\gamma)^i x_i },
\]
so the upper bound is trivial. The lower bound follows in the same way
as the lower bound for $b_j^0(x)$ once we write
\[
d_j^0(x) = x_j \frac{x_m + \sum_{i\neq j,m} (1+\gamma)^{i-m} x_i}{N
+ \sum_i [(1+\gamma)^{i-m} - 1 ] x_i} \geq x_j \frac{N - \sum
_{i\neq m} x_i}{N + \sum_i [(1+\gamma)^{i-m} - 1 ] x_i}.
\]
\upqed\end{pf}

\begin{table}[b]
\tablewidth=10,7cm
\caption{Rates for the comparison branching processes, $j \geq1$}
\label{table:compbp}
\begin{tabular*}{10,7cm}{@{\extracolsep{4in minus 4in}}lcc@{}}
\hline
& $\bolds{Z_{j,u}^\mu(s)}$ & $\bolds{Z_{j,\ell}^\mu(s)}$\\
\hline
Birth rate & $b_{j,u}^\mu\equiv\frac{(1+\rho)(1+\gamma
)^{j-m}}{1-h_\mu
}$ & $b_{j,\ell}^\mu\equiv\frac{(1+\rho-\mu^\eta)(1+\gamma
)^{j-m}}{1+g_\mu}$
\\
Death rate & $d_{j,u}^\mu\equiv\frac{1-M\mu^\eta}{1+g_\mu}$ &
$d_{j,\ell}^\mu\equiv1+\mu$
\\
Immigration rate & $\mu Z_{j-1,u}^\mu(t)$ & $\mu Z_{j-1,\ell}^\mu(t)$
\\
\hline
\end{tabular*}
\end{table}

We now describe the bounding processes. Let $0 < t_1 < t_2 < \Delta_0$,
\[
\eta= \eta(t_2) = \frac{\lambda_1-\rho}{4\gamma}(\Delta_0-t_2).
\]
The reason for this choice of $\eta$ is that
\[
y_j(t) \leq(\alpha+t \rho/\gamma) -4\eta
\]
for all $t \leq t_2$, $j \geq1$. For our bounding processes, we set
$Z_{0,u}(t) \equiv N^\mu(t)$, $Z_{0,\ell}^\mu(t) \equiv(1-k\mu
^\eta) N^\mu(t) $ and let $Z_{j,a}^\mu$, $ 1 \leq j \leq k$,
$a=u,\ell$, be (birth and death) branching processes with rates given
in Table \ref{table:compbp}, taking $m=0$, $M=k$. Note that the birth
and death rates are per particle. The extra factor $\mu$ in the
definition of $d_{j,\ell}^\mu$ takes care of deaths due to mutations.
We also set $\lambda_{j,a}^\mu\equiv b_{j,a}^\mu- d_{j,a}^\mu$ to
be the growth rates of $Z_{j,a}^\mu$, $a=u,\ell,$ so that we have
$\lambda_{j,a}^\mu\to\lambda_j$ as $\mu\to0$ for $j \geq1$,
$a=u,\ell$. If we use the convention that $\lambda_{0,a}^\mu= \rho$
for $a=u,\ell$, this also holds for $j=0$.

For the next result, we use the notation $Z_a^\mu(t) = (Z_{0,a}^\mu
(t),Z_{1,a}^\mu(t), \ldots, Z_{k,a}^\mu(t),\break 0,\ldots)$, for $a=u,\ell$.

\begin{lemma} \label{coupling}
There exists a coupling of $X^\mu(t)$ with $Z_{a}^\mu(t)$, $a=u,\ell
,$ such that
\[
Z_{j,\ell}^\mu(t) \leq X_j^\mu(t) \leq Z_{j,u}^\mu(t)
\]
for all $t \leq( L \sigma/\gamma)\wedge T_{k+1}^\mu$, $j \leq k$, where
\begin{eqnarray*}
\sigma&=& \inf\{t \geq0\dvtx   Y_i^\mu(t)>\alpha+t \rho/\gamma
-2\eta
\\
&&{}\qquad\hspace*{-6pt}\mbox{ for some } i \geq1   \mbox{ or } |F^\mu(t) - (\alpha+ t
\rho/\gamma)| > \eta\}.
\end{eqnarray*}
\end{lemma}

\begin{pf}
For $t \leq T_{k+1}^\mu$, we have $X_j(t) = 0$ if $j > k$.
Furthermore, if $t \leq\sigma L/\gamma$,
\[
\frac{N^\mu(Lt/\gamma)}{(1/\mu)^\alpha e^{\rho t}} \leq(1/\mu
)^\eta
\]
so that
\[
\frac{X_j^\mu(Lt/\gamma)}{N^\mu(Lt/\gamma)} \leq\frac{\mu^2\eta
}{(1/\mu)^\eta}= \mu^\eta
\]
for all $j \geq1$ and hence we have the bounds on birth and death
rates given in Lemma \ref{bdbnds} with $m=0$ and $M=k$. The processes
can therefore be coupled in an elementary way by matching birth, deaths
and immigrations in the appropriate manner.
\end{pf}

The result which we will dedicate most of the remainder of this section
to proving is the following.

\begin{lemma} \label{propcompbp}
Let $ 0 < t_1 < t_2 < \delta_0$. For $a=u,\ell$ and $j \leq k$, we have
\[
(1/L)\log^+ Z_{j,a}^\mu(Lt/\gamma) \to y_j(t)
\]
in probability uniformly on $[t_1,t_2]$.
\end{lemma}

Because $y(t) \leq(\alpha+\rho t) -4\eta$ for all $t \leq t_2$,
Lemma \ref{propcompbp} implies that
\[
P\bigl( (1/L) \log^+ Z_a^\mu(Lt/\gamma) \leq(\alpha+\rho t) -2\eta,
\forall  t \leq t_2,   a=u,\ell\bigr) \to1
\]
as $\mu\to0$. This and Lemma \ref{poplim} imply that $P(\sigma>
t_2) \to1 $ as $\mu\to0$ and therefore Proposition \ref{startup}
follows from Lemmas \ref{propcompbp},  \ref{coupling} and
\ref{nobettertype}.

To prove Lemma \ref{propcompbp}, we begin by defining another level of
upper and lower bounds, $\hat{Z}_{j,a}^\mu$, in which immigrations
occur at deterministic rates. More specifically, for $a=u,\ell$, we
define $\hat{Z}_{j,a}^\mu(t)$ as a branching process with the same
initial population and birth and death rates as $Z_{j,a}^\mu(t)$, but
with immigrations at rate $\mu I_{j,a}^\mu(t),$ where
\[
I_{j,u}^\mu(t) \equiv E \hat{Z}_{j-1,u}^\mu(t) + e^{\lambda
_{j-1,u}^\mu t} (1/\mu)^{2(\alpha-(j-1))/3}
\]
and
\[
I_{j,\ell}^\mu(t) \equiv E \hat{Z}_{j-1,\ell}^\mu(t) - e^{\lambda
_{j-1,\ell}^\mu t} (1/\mu)^{2(\alpha-(j-1))/3}.
\]
We will use the convention that $I_{0,a}^\mu(t) \equiv0$ for all $t$.
Note that
%
\begin{equation}\label{mexpress}
E(e^{-\lambda_{j,a}^\mu t}\hat{Z}_{j,a}^\mu(t)) = \mu\int_0^t
e^{-\lambda_{j,a}^\mu s} I_{j,a}^\mu(s) \, ds
\end{equation}
for all $j \geq1$ and $a=u,\ell$, an expression which will be used
often throughout the remainder of this section.

\begin{lemma} \label{Zjmart}
For $j \geq0$ and $a=u,\ell$,
\[
M_{j,a}^\mu(t) \equiv e^{-\lambda_{j,a}^\mu t}\hat{Z}_{j,a}^\mu(t)
- E(e^{-\lambda_{j,a}^\mu t}\hat{Z}_{j,a}^\mu(t))
\]
is a martingale with respect to the filtration
\[
\mathcal{G}_{a,t} \equiv\sigma\{\hat{Z}_{i,a}^\mu(s)\dvtx 0 \le i \le j,
  s \le t\}.
\]
\end{lemma}

\begin{pf}
We prove the result for $a=u$, the proof for $a=\ell$ being similar,
and drop the subscripts $u$ from all quantities for the remainder of
the proof. It is easy to see that
\[
E\bigl( \hat{Z}_j^\mu(t+h) | \mathcal{G}_t \bigr) = e^{\lambda_j^\mu h} \hat
{Z}_j^\mu(t)
+ E \biggl(\mu\int_t^{t+h} e^{\lambda_j^\mu(t+h-s)}
I_{j}^\mu(s) \, ds \big|\mathcal{G}_t \biggr)
\]
and multiplying by $e^{-\lambda_j^\mu(t+h)}$ gives
\[
E \biggl( e^{-\lambda_j^\mu(t+h)} \hat{Z}_j^\mu(t+h) -
e^{-\lambda_j^\mu t} \hat{Z}_j^\mu(t) -
\mu\int_t^{t+h} e^{-\lambda_j^\mu s}I_{j}^\mu(s) \,  ds
\big|\mathcal{G}_t \biggr) = 0.
\]
Since \eqref{mexpress} implies that
\begin{eqnarray*}
&&M_j^\mu(t+h) - M_j^\mu(t)
\\
&&\qquad = e^{-\lambda_j^\mu(t+h)} \hat{Z}_j^\mu
(t+h) - e^{-\lambda_j^\mu(t)} \hat{Z}_j^\mu(t)
- \mu\int_t^{t+h}e^{-\lambda_j^\mu s} I_{j}^\mu(t)   \,ds
\end{eqnarray*}
for $j \geq1$ and the same equality clearly holds for $j = 0$ as well,
the desired result, $E(M_j^\mu(t+h)-M_j^\mu(t)|\mathcal{G}_t) = 0$, follows.
\end{pf}

\begin{lemma} \label{Zmgbnd}
For all $a=u,\ell$, $T > 0$ and $\mu$ sufficiently small, we have
\begin{eqnarray*}
P \Bigl( \sup_{t\le T} |M_{j,a}^\mu(t)| > (1/\mu)^{2(\alpha-j)/3}
\Bigr) \le C\mu^{(\alpha-j)/3} \bigl[1 + \mu^{(\alpha-j+1)/3}\bigr].
\end{eqnarray*}
In particular, for all $j \leq k,$
\begin{eqnarray*}
P\bigl( |\hat{Z}^\mu_{j,a}(t) - E \hat{Z}^{\mu}_{j,a}(t)| > e^{\lambda
_{j,a}^\mu t} (1/\mu)^{2(\alpha-j)/3},   \forall  t \leq T\bigr) \to0
\end{eqnarray*}
as $\mu\to0$
\end{lemma}

\begin{pf}
The second part of the result follows directly from the first, along
with the definition of $M_{j,a}^\mu(t)$. To obtain the first part, we
suppose for the remainder of the proof that $u=a$ and drop the
subscript $u$. The proof for $a = \ell$ is similar. We will also
assume that $j \geq1$ and leave the (simpler) $j=0$ case to the reader.

We proceed by calculating the variance of $e^{-\lambda_j^\mu t}\hat
{Z}_j^\mu(t)$ and then using the $L^2$ maximum inequality to bound the
second moment of $M_j^\mu(t)$ uniformly on $[0,T]$. To begin, we claim
that provided we choose $\mu$ small enough so that $\lambda_i^\mu>
\lambda_{i-1}^\mu$ for all $1 \leq i \leq j$, we have
%
\begin{equation}\label{meanbnd}
\qquad g(t) (1/\mu)^{(\alpha-j)} \leq E(e^{-\lambda_j^\mu t}\hat{Z}_j^\mu
(t)) \leq C \bigl[(1/\mu)^{(\alpha-j)} + \mu^{1/3}(1/\mu)^{2(\alpha-j)/3}\bigr],\hspace*{-8pt}
\end{equation}
where $g(t)$ is continuous on $[0,\infty)$ and positive on $(0,\infty
)$. To see this, we note that $E \hat{Z}_0(t) = E N^\mu(t) = (1/\mu
)^\alpha e^{\rho t}$, so the result clearly holds for $j=1$ by \eqref
{mexpress} and the general case follows by induction on $j$. Now,
\begin{eqnarray} \label{varexpress}
\quad &&\frac{d}{dt}E(e^{-\lambda_j^\mu t}\hat{Z}_j^\mu(t))^2\nonumber
\\
&&\qquad = -2\lambda_j^\mu E(e^{-\lambda_j^\mu t}\hat{Z}_j^\mu(t))^2 +
e^{-2 \lambda_j^\mu t}E\bigl[b_j^\mu\hat{Z}_j^\mu(t)\bigl(2\hat{Z}_j^\mu
(t)+1\bigr)\bigr] \nonumber
\\[-8pt]\\[-8pt]
&& {}\quad\qquad-d_j^\mu e^{-2 \lambda_j^\mu t}E\bigl[\hat{Z}_j^\mu(t)\bigl(2\hat
{Z}_j^\mu(t)-1\bigr)\bigr]+ \mu I_j^\mu(t) e^{-2 \lambda_j^\mu t}E[2\hat
{Z}_j^\mu(t)+1]\nonumber
\\
&&\qquad = (b_j^\mu+ d_j^\mu)e^{-2\lambda_j^\mu t}E\hat{Z}_j^\mu(t) + \mu
I_j^\mu(t) e^{-2 \lambda_j^\mu t} + 2 \mu I_j^\mu(t) e^{-2 \lambda
_j^\mu t}E \hat{Z}_j^\mu(t).\nonumber
\end{eqnarray}
Equation (\ref{mexpress}) implies that
\begin{eqnarray*}
&& \int_0^t 2\mu I_j^\mu(s) e^{-2\lambda_j^\mu s}E\hat{Z}_j^\mu(s) \,ds
  \\
&& \qquad= 2 \int_0^t \mu I_j^\mu(s)e^{-\lambda_j^\mu s} \int_0^s
\mu I_j^\mu(r) e^{-\lambda_j^\mu r} \, dr  \, ds
= [E (e^{-\lambda_j^\mu t} \hat{Z}_j^\mu(t))]^2
\end{eqnarray*}
so that integrating both sides of \eqref{varexpress} and applying
\eqref{meanbnd} yields
\begin{eqnarray}\label{varZjt}
\quad \operatorname{var}(e^{-\lambda_j^\mu t}\hat{Z}_j^\mu(t)) & \le&
(b_j^\mu+ d_j^\mu
) \int_0^t e^{-2\lambda_j^\mu s} E\hat{Z}_j^\mu(s)  \, ds
+ \int_0^t \mu I_j^\mu(s) e^{-2\lambda_j^\mu s} \,  ds
\nonumber
\\[-8pt]\\[-8pt]
& \le& C \bigl[(1/\mu)^{(\alpha-j)} + \mu^{1/3}(1/\mu)^{2(\alpha-j)/3}\bigr].\nonumber
\end{eqnarray}
By Lemma \ref{Zjmart}, $M_j^\mu$ is a martingale with respect to
$\mathcal{G}_t$ and so the $L^2$ maximum inequality implies that
\[
E \Bigl( \sup_{t\le T} (M_j^\mu(t))^2 \Bigr) \le4 E(M_j^\mu
(T))^2 = 4 \operatorname{var}(e^{-\lambda^\mu_j t_2}\hat{Z}_j^\mu(T)),
\]
the second equality following from the definition of $M_j^\mu$.
Applying Chebyshev's inequality and \eqref{varZjt} then yields
\[
P \Bigl( \sup_{t\le t_2} |M_j^\mu(t)| > (1/\mu)^{2(\alpha-j)/3}
\Bigr) \le C \mu^{(\alpha-j)/3} \bigl[1 + \mu^{(\alpha-j+1)/3}\bigr],
\]
completing the proof.
\end{pf}

\begin{corollary} \label{Zmgbndcor}
For $a=u,\ell$, there exists a coupling of the process $Z_a^\mu$ with
$\hat{Z}_a^\mu$ such that
\[
\hat{Z}_{j,\ell}^\mu(t) \leq Z_{j,\ell}^\mu(t) \leq Z_{j,u}^\mu
(t) \leq\hat{Z}_{j,u}^\mu(t)
\]
for all $t \leq\hat{\sigma,}$ where
\[
\hat{\sigma} \equiv\inf\{t \geq0\dvtx \hat{Z}_{j-1,u}^\mu(t) >
I_{j,u}^\mu(t) \mbox{ or } \hat{Z}_{j-1,\ell}^\mu(t) < I_{j,\ell
}^\mu(t) \mbox{ for some } j \geq1\}.
\]
Furthermore, $P(\hat{\sigma} \leq t_2) \to0$ as $\mu\to0$.
\end{corollary}

\begin{pf}
Arguing inductively, we can see that the immigration rates for type
$j$'s in $\hat{Z}_{\ell}^\mu, Z_{\ell}^\mu, Z_{u}^\mu$ and $\hat
{Z}_{u}^\mu$, respectively, satisfy
\[
\mu I_{j,\ell}^\mu(t) \leq\mu Z_{j-1,\ell}^\mu(t) \leq\mu
Z_{j-1,u}^\mu(t) \leq\mu I_{j,u}^\mu(t)
\]
for $t \leq\hat{\sigma}$. Therefore, we define a coupling for the
two processes by coupling births, deaths and immigrations. The fact
that $P(\hat{\sigma} \leq t_2) \to0$ follows from Lemma \ref{Zmgbnd}.
\end{pf}

Define $\hat{Y}_{j,a}^\mu(t) \equiv(1/L)\log^+ \hat{Z}_{j,a}^\mu
(Lt/\gamma)$. Lemma \ref{propcompbp} follows from Corollary \ref
{Zmgbndcor} along with our next result.

\begin{lemma} \label{Zmgcor}
Let $j \leq k$, $a=u,\ell$. Then, $\hat{Y}_{j,a}^\mu(t) \to y_j(t)$
in probability uniformly on $[t_1,t_2]$.
\end{lemma}

\begin{pf}
Again, we only prove the result for $u=a$ and drop the $u$ subscript.
Let $j \leq k$ and write
%
\begin{eqnarray} \label{Yhatlim}
&&\hat{Y}_j^\mu(t) - y_j(t)\nonumber
\\
&&\qquad = \bigl((1/L)\log^+ [e^{-\lambda_j^\mu Lt/\gamma} \hat{Z}_j^\mu
(Lt/\gamma)]\nonumber
\\[-8pt]\\[-8pt]
&&\hspace*{1,5pt}{}\qquad\quad - (1/L)\log^+ E [e^{-\lambda_j^\mu Lt/\gamma} \hat
{Z}_j^\mu
(Lt/\gamma)] \bigr)\nonumber
\\
&&{} \hspace*{-2,2pt}\qquad\quad+ \bigl((1/L)\log^+ E [e^{-\lambda_j^\mu Lt /\gamma} \hat
{Z}_j^\mu(Lt/\gamma)] - (\alpha- j)\bigr)
+ (\lambda_j^\mu- \lambda_j)t/\gamma.\nonumber
\end{eqnarray}
By Lemma \ref{Zmgbnd},
\begin{eqnarray*}
P \bigl( |\hat{Z}_{j,u}^\mu(t) - E \hat{Z}_{j,u}^\mu(t)| \leq
e^{\lambda_j^\mu t} (1/\mu)^{2(\alpha-j)/3},   \forall  t \leq
t_2 \bigr) \to1,
\end{eqnarray*}
and on the set where
\[
|\hat{Z}_{j,u}^\mu(t) - E \hat{Z}_{j,u}^\mu(t)| \leq e^{\lambda
_j^\mu s} (1/\mu)^{2(\alpha-j)/3} \qquad\forall  t \leq t_2,\qquad
\]
we have
\begin{eqnarray*}
&&(1/L)\log^+ [e^{-\lambda_j^\mu sL/\gamma} \hat{Z}_j^\mu(Lt/\gamma
)] -
(1/L)\log^+ E [e^{-\lambda_j^\mu Lt/\gamma} \hat{Z}_j^\mu
(Lt/\gamma)]
\\
&&\qquad = \frac{1}{L} \log\biggl(1 + \frac{\hat{Z}_j^\mu(Lt/\gamma) - E
\hat{Z}_j^\mu(Lt/\gamma)}{E \hat{Z}_j^\mu(Lt/\gamma)}\biggr)
\\
&&\qquad \leq\frac{C}{L} \frac{|\hat{Z}_j^\mu(Lt/\gamma) - E \hat
{Z}_j^\mu
(Lt/\gamma)|}{E \hat{Z}_j^\mu(Lt/\gamma)}
\leq\frac{C}{L} \frac{(1/\mu)^{2(\alpha-j)/3}}{(1/\mu)^{\alpha
-j}} \to0
\end{eqnarray*}
uniformly on $[t_1,t_2]$ as $\mu\to0$, the last inequality following
from \eqref{meanbnd} and the fact that $g(t)$ is bounded away from 0
on $[t_1,t_2]$. Therefore, the absolute value of the first term on the
right of \eqref{Yhatlim} goes to zero uniformly on $[t_1,t_2]$. It is
clear from \eqref{meanbnd} that the second term goes to 0 as well and
since $\lambda_j^\mu\to\lambda_j$ as $\mu\to0$, the result follows.
\end{pf}

\section{Inductive step} \label{sec-inductpf}

In this section we prove Propositions \ref{interior}--\ref{case2}. We
shall assume throughout that the assumptions from Section \ref
{pf-ideas} hold and begin with the proof of Proposition \ref
{interior}. The reader should refer to the statement of that result for
the notation used throughout this section.

\subsection{Interior convergence}

Let $\varepsilon>0$, set $a_j(t) \equiv\alpha+ t\rho/\gamma$ for
$j \neq k$,
$a_k(t) \equiv1$ and choose $\eta= \eta(\varepsilon) >0$ so that: (i)
$y_j(t) < a_j(t) - 2\eta$ for all $t \leq\Delta_\varepsilon$, $j
\neq m;$
(ii) $y_{j-1}(t)-y_j(t) < 1-2\eta$ for all $t \leq\Delta_\varepsilon
$, $j
\geq0$. Given $\zeta>0$, we define the stopping times
\begin{eqnarray*}
\sigma_0(j) &\equiv&\gamma T_{j}^\mu/L,
\\
\sigma_1(j) &\equiv&\inf\{ t \geq0\dvtx Y^\mu_j(t) \ge a_j(t) - \eta\}
,\\
\sigma_1 &\equiv&\inf_{j \neq m} \sigma_1(j),
\\
\sigma_1' &\equiv&\inf_{j < m} \sigma_1(j),
\\
\sigma_2(j) &\equiv&\inf\{t \geq0\dvtx Y^\mu_{i-1}(t)-Y^\mu_i(t) \ge1
- \eta,   \mbox{ for some } 1 \leq i \leq j \},
\\
\sigma_3(j,\zeta) &\equiv&\inf\{ t \geq0\dvtx Y^\mu_j(t) \le\zeta\}.
\end{eqnarray*}
For the remainder of this section, set $\sigma_0 = \sigma_0(k+1)$ and
$\sigma_2 = \sigma_2(k)$. We shall prove convergence of $Y_j^\mu(t)$
up to time $\sigma(j,\zeta) \equiv\sigma_0 \wedge\sigma_1 \wedge
\sigma_2 \wedge\sigma_3(j,\zeta)$. For types $j \leq k$, this will
essentially amount to controlling the infinitesimal variance of
$Y_j^\mu$ (Lemma \ref{martbnd}) and then showing that the
infinitesimal mean converges to the appropriate limit (Lemma \ref
{infmeanlim}), while for types $j >k$ we will simply show that they are
unlikely to be born before time $L \Delta_\varepsilon/\gamma$, that
is, $\sigma
_0 > \Delta_\varepsilon$ with high probability (this follows from
Lemma~\ref
{nonewguys}). We then complete the proof of Proposition \ref{interior}
by using the structure of the limit $y_j(t)$ to extend convergence up
to time $\Delta_\varepsilon$, as required. If $y_j(t)$ is bounded
away from
0, then this is easy since our choice of $\eta$ implies that $\sigma
_1, \sigma_2$ are unlikely to occur before time $\Delta_\varepsilon$
and if
$y_j(t)$ is not bounded away from~0 (which can only happen if $j < m$),
we will define a stopping time $\sigma'$ such that $Y_j^\mu(t \wedge
\sigma')$ is a supermartingale to conclude that once the $j$'s drop
below a certain level, they will never climb up again.


The first step is to calculate the infinitesimal mean and variance.
Writing $y_j = (1/L) \log(x_j)$, $y = (y_0,y_1,\ldots)$ and $N = \sum
e^{Ly_i}$, noting the time rescaling and using the fact that the change
in $y_j$ when $x_j$ jumps to $x_j \pm1$ is $(1/L)\log(1\pm
x_j^{-1})$, we can write the infinitesimal mean of $Y_j^\mu(t)$ as
$B_j(y) = B_{j,r}(y) + B_{j,m}(y),$ where
\begin{eqnarray*}
B_{j,r}(y) & = &\gamma^{-1} \frac{[(1+\rho)N-e^{Ly_j}](1+\gamma)^j
e^{Ly_j}}{\sum_{i \geq0} (1+\gamma)^i e^{Ly_i}} \log(1+e^{-Ly_j})
\\
&& {}+ \gamma^{-1}\frac{\sum_{i \neq j} (1+\gamma)^i
e^{Ly_i}}{\sum_{i
\geq0} (1+\gamma)^i e^{Ly_i}} e^{Ly_j} \log(1-e^{-Ly_j}),
\\
B_{j,\mu} & =& \mu\gamma^{-1} e^{Ly_{j-i}}\log(1+e^{-Ly_j}) + \mu
\gamma
^{-1} e^{Ly_{j}} \log(1-e^{-Ly_j}).
\end{eqnarray*}
In words, $B_{j,r}(y)$ is the rate of change due to death and
subsequent replacement, while $B_{j,\mu}(y)$ is the rate of change due
to mutations. Similarly, the infinitesimal variance is
\begin{eqnarray*}
A_{j}(y) &= & (1/L) \biggl[ \gamma^{-1}\frac{[(1+\rho
)N-e^{Ly_j}](1+\gamma
)^j e^{Ly_j}}{\sum_{i \geq0} (1+\gamma)^i e^{Ly_i}} \bigl(\log
(1+e^{-Ly_j})\bigr)^2
 \\
&&{} \qquad\hspace*{8pt}+ \gamma^{-1}\frac{\sum_{i \neq j} (1+\gamma)^i
e^{Ly_i}}{\sum_{i
\geq0} (1+\gamma)^i e^{Ly_i}} e^{Ly_j} \bigl(\log(1-e^{-Ly_j})\bigr)^2
\\
&&{} \qquad\hspace*{8pt}+ \mu\gamma^{-1} e^{Ly_{j-i}}\bigl(\log(1+e^{-Ly_j})\bigr)^2 + \mu
\gamma
^{-1} e^{Ly_{j}} \bigl(\log(1-e^{-Ly_j})\bigr)^2 \biggr].
\end{eqnarray*}
Introducing $f_1(x) \equiv x \log(1+x^{-1})$, $f_2(x) \equiv x \log
(1-x^{-1})$,
\begin{eqnarray*}
g_{j,1}(y) \equiv\frac{(1+\gamma)^j}{\gamma} \frac{[(1+\rho
)N-e^{Ly_j}]}{\sum_{i \geq0} (1+\gamma)^i e^{Ly_i}}
\end{eqnarray*}
and
\begin{eqnarray*}
g_{j,2}(y) \equiv\frac{1}{\gamma} \frac{\sum
_{i\neq
j} (1+\gamma)^i e^{Ly_i}}{\sum_{i \geq0} (1+\gamma)^i e^{Ly_i}},
\end{eqnarray*}
we can write
\begin{eqnarray*}
A_j(y) &=& (1/L)  [e^{-Ly_j} g_{j,1}(y) f_1^2(e^{Ly_j})+
e^{-Ly_j} g_{j,2}(y) f_2^2(e^{Ly_j})
\\
&&{}\hspace*{16pt} \quad+ \mu\gamma^{-1} f_1^2(e^{Ly_j})e^{Ly_{j-i}-2Ly_j} + \mu
\gamma^{-1}
f_2^2(e^{Ly_j})e^{-Ly_j}].
\end{eqnarray*}
Since $g_{j,1}(x) \leq(1+\rho)(1+\gamma)^j/\gamma$, $g_{j,2}(x)
\leq1/\gamma$
and $f_2(x) \leq f_1(x) \leq1$ for all $x \in[0,\infty)$, we obtain
the bound
%
\begin{equation}\label{infvarbnd}
A_j(y) \leq(C/L)\bigl((1+\mu)e^{-L y_j} + \mu e^{Ly_{j-i}-2Ly_j}\bigr).
\end{equation}
Define
\[
M_j(t) = Y_j^\mu(t)-Y_j^\mu(0) - \int_0^t B_j(Y^\mu(s))\,ds.
\]

\begin{lemma} \label{martbnd}
For any $\xi>0$ and $j \geq1$,
\[
P\Bigl(\sup_{t \leq\sigma_2(j)} |M_j(t)|> \xi\Bigr) \to0.
\]
\end{lemma}

\begin{pf}
Since $Y_{j-1}^\mu(t) - Y_j^\mu(t) < 1-\eta$ for $t \le\sigma_2$,
we have
\[
e^{-L Y_j^\mu(t)} \leq1
\quad\mbox{and}\quad
\mu e^{LY_{j-i}^\mu(t)- 2LY_j^\mu(t)} \leq C \mu^{\eta}
\]
and therefore the result follows from \eqref{infvarbnd} and Corollary
2.8 in \cite{TK}.
\end{pf}

Our next step is to show that the infinitesimal means converge to the
appropriate limit. The key to the proof is that $b_j^\mu(Y_j^\mu
(t))\to\lambda_{j-m}\gamma$ for all $t \leq\sigma(j,\zeta)$, but we
write out the details carefully because we will need (I)--(IV) from the
proof several times in what follows.

\begin{lemma} \label{infmeanlim} If $j \neq m$ and $\zeta, \xi>0,$
then as $\mu\to0$,
\[
P\biggl(\sup_{t \leq\sigma(j,\zeta)}\bigg|\int_0^t B_{j}(Y_j^\mu
(s)) \, ds - \lambda_{j-m}t/\gamma\bigg|> \xi\biggr) \to0.
\]
\end{lemma}

\begin{pf}
Using the definition of $f_i$, $g_{j,i}$, $i=1,2$, we write
\[
B_{j,r}(y) = f_1(e^{Ly_j}) g_{j,1}(y) + f_2(e^{Ly_j}) g_{j,2}(y)
\]
and
\[
B_{j,\mu}(y) = \mu\gamma^{-1}\bigl[ f_1(e^{Ly_j}) e^{L(y_{j-1}-y_j)} +
f_2(e^{Ly_j})\bigr].
\]
We will complete the proof by proving the following four facts:
\begin{enumerate}[(III)]
\item[(I)] for any $\zeta>0$, $f_1(e^{LY_j^\mu(t)}) \to1$ and
$f_2(e^{LY_j^\mu(t)}) \to-1$ in probability uniformly on $[0,\sigma
_3(j,\zeta)]$;
\item[(II)] $g_{j,2}(Y^\mu_j(t)) \to1/\gamma$ in probability uniformly
on $[0,\sigma_1(j)];$
\item[(III)] for any $\zeta>0$,
\[
P \Bigl( \sup_{t \leq\sigma_1'} g_{j,1}(Y^\mu_j(t)) > (1+\rho
)(1+\gamma)^{j-m}/\gamma+ \zeta\Bigr) \to0
\]
and, furthermore, $g_{j,1}(Y^\mu_j(t)) \to(1+\rho)(1+\gamma
)^{j-m}/\gamma$
in probability uniformly on $[0,\sigma_0 \wedge\sigma_1];$
\item[(IV)] $B_{j,\mu}(Y^\mu(t)) \to0$ in probability uniformly on
$[0,\sigma_2(j)].$
\end{enumerate}
(I) follows immediately since $f_1(x) \to1$, $f_2(x) \to-1$ as $x \to
\infty$ and $Y_j^\mu(t) \geq\zeta$ on $[0,\sigma_3(j,\zeta)]$. To
prove (II), write
\begin{eqnarray*}
g_{j,2}(y) = \frac{1}{\gamma}\biggl(1 - \frac{(1+\gamma)^j
e^{Ly_j}}{\sum
_{i \geq0} (1+\gamma)^i e^{Ly_i}}\biggr)
\end{eqnarray*}
and note that if $\sum e^{Ly_i} = N$ and $y_i \leq a$, then
\[
0 \leq\frac{(1+\gamma)^j e^{Ly_j}}{\sum_{i \geq0} (1+\gamma)^i e^{Ly_i}}
\leq(1+\gamma)^j e^{Ly_j}/N \leq C (1/\mu)^{a}/N.
\]
Now, Lemma \ref{poplim} and the assumption that $F^\mu(0) \to\alpha
$ imply that $N^\mu(Lt/\gamma) \geq(1/\mu)^{\alpha+t\rho/\gamma
-\eta
/2}$ for all $t \leq\Delta_0$ with high probability so that since
$\sum e^{L Y_i^\mu(t)} = N^\mu(t)$ and $Y_j^\mu(t) \leq\alpha+ t
\rho/\gamma-\eta,$ if $t \leq\sigma_1(j)$, (II) follows. For (III) we
note that if $\sum e^{Ly_i} = N$, then using the definition of $\gamma
_j =
(1+\gamma)^j -1$, we have
\begin{eqnarray*}
g_{j,1}(y) = \frac{(1+\gamma)^{j-m}}{\gamma}\biggl(\frac{1+\rho
-e^{Ly_j}/N}{1 + \sum_{i \neq m} \gamma_{i-m} e^{Ly_i}/N}\biggr).
\end{eqnarray*}
The first part of (III) then follows from the fact that
\[
\sum_{i \neq m} \gamma_{i-m} e^{Ly_i}/N \geq\sum_{i < m} \gamma_{i-m}
e^{Ly_i}/N \geq\gamma_{-m}(1/\mu)^a/N
\]
if $y_i \leq a$ for all $i < m,$ while the second part follows from the
fact that we also have
\[
\biggl(\frac{1+\rho-e^{Ly_j}/N}{1 + \sum_{i \neq m} \gamma_{i-m}
e^{Ly_i}/N}\biggr) \geq\frac{1+\rho-(1/\mu)^a/N}{1 + (k-m)\gamma
_{k-m}(1/\mu)^a/N}
\]
if $y_j \leq a$ for all $j \leq k$ and $y_j = 0$ for all $j >k$.
Finally, to prove (IV), we use the bound
\[
B_{j,\mu}(y) \leq C \mu\bigl[f_1(e^{Ly_j}) e^{L(y_{j-1}-y_j)} + f_2(e^{Ly_j})\bigr]
\]
so that since $f_2(x), f_1(x) \leq1$ for all $x \geq0$, the result
follows from the fact that
\[
\mu e^{L(y_{j-1}-y_j)} \leq\mu^{\eta}
\]
if $y_{j-1}-y_j < 1-\eta$.
\end{pf}

%
%
%

\begin{pf*}{Proof of Proposition \ref{interior}}
Lemma \ref{nonewguys} implies that if we have the result for $j \leq
k$, then $P(\sigma_0 \leq\Delta_\varepsilon) \to0,$ so it suffices to
prove the result for $j \leq k$. Lemmas~\ref{martbnd} and \ref
{infmeanlim} and the assumption that $Y_j^\mu(0) \to y_j^0$ in
probability together imply that $Y_j(t) \to y_j(t)$ in probability
uniformly on $[0,\sigma(j,\zeta)]$ for any $\zeta>0$. Note that
since $y_i(t) \leq a_i(t) - 2\eta$ for all $t \leq\Delta_\varepsilon
$, $i
\neq m,$ by our choice of $\eta$, we have $P(\sigma_1 \leq\Delta
_\varepsilon\wedge\sigma_0 \wedge\sigma_2) \to0$. Furthermore,
$y_k(t) <
1-\varepsilon$ for all $t \leq\Delta_\varepsilon$ and hence
$P(\sigma_0 \leq
\Delta_\varepsilon\wedge\sigma_2) \to0$ by Lemma \ref{nonewguys}.
Therefore, we obtain uniform convergence on $[0,\Delta_\varepsilon
\wedge
\sigma_2 \wedge\sigma_3(j,\zeta)]$ for any $\zeta>0$. We will show
that the convergence is uniform on $[0,\Delta_\varepsilon\wedge
\sigma_2]$.
Proposition \ref{interior} follows since $y_{j-1}(t) - y_j(t) <
1-2\eta$ for all $t \leq\Delta_\varepsilon$, $j \geq0$.

Suppose first that $j \geq m$. There then exists $\zeta>0$ such that
$y_j(t) \geq\zeta$ for all $t \leq\Delta_\varepsilon$ and hence we obtain
uniform convergence on $[0, \Delta_\varepsilon\wedge\sigma_2(j)]$.
The same
argument applies if $j < m$ and $y_j(t)$ is bounded away from 0 on
$[0,\Delta_\varepsilon]$. If $j < m$ and $y_j(t)$ is not bounded away
from 0
on $[0,\Delta_\varepsilon]$, set
\[
\sigma' = \sigma_0\wedge\sigma_1 \wedge\sigma_2,
\]
let $\zeta, \xi> 0$ be small and choose a time $t_0 < \Delta
_\varepsilon$
such that $y_j(t_0) = \zeta\xi/4$. If no such time exists [i.e.,
$y_j(t)= 0$ for all $t$], set $t_0 = 0$.
Then, $y_j(s) \leq\zeta\xi/4$ for all $t_0\leq s \leq\Delta
_\varepsilon$ and
\begin{eqnarray} \label{smallguybnd}
&&P\Bigl(\sup_{t \leq\Delta_\varepsilon\wedge\sigma'} |Y_j(t)
 -y_j(t) | > \zeta\Bigr) \nonumber
\\
&&\qquad \leq P\Bigl( \sup_{t \leq t_0\wedge\sigma'} |Y_j(t)-y_j(t)| >
\zeta\Bigr) + P\Bigl( \sup_{t_0 \leq t \leq\Delta_\varepsilon
}Y_j(t\wedge\sigma') > \zeta(1+\xi/4)\Bigr) \nonumber
\\[-8pt]\\[-8pt]
&&\qquad\leq P\Bigl( \sup_{t \leq t_0\wedge\sigma'} |Y_j(t)-y_j(t)| >
\zeta\Bigr) + P\bigl(Y_j(t_0\wedge\sigma') > \zeta\xi/2\bigr)
\nonumber
\\
&&{}\qquad \quad+ P\Bigl(\sup_{t_0 \leq t \leq\Delta_\varepsilon}
Y_j^\mu
(t\wedge\sigma') > \zeta(1+\xi/4) \big|Y^\mu_j(t_0\wedge\sigma
') \le\zeta\xi/2 \Bigr).\nonumber
\end{eqnarray}
The argument in the last paragraph implies that $Y_j^\mu(t) \to
y_j(t)$ uniformly on $[0,t_0 \wedge\sigma']$ and hence the first and
second terms on the right-hand side of \eqref{smallguybnd} are each $<
\xi/4$ for all $\mu$ sufficiently small. To control the third term,
we note that (II)--(IV) from the proof of Lemma \ref{infmeanlim},
along with the bounds $f_1(x) \leq1$, $f_2(x) \leq-1$ for all $x \geq
0,$ imply that if $\mu$ is sufficiently small, then $B_j(Y^\mu(t))
\leq0$ for all $t \leq\sigma'$ with high probability so that $Y^\mu
(t\wedge\sigma')$ is a supermartingale. Therefore,
\begin{eqnarray*}
P\Bigl(  \max_{t_0\leq t \leq\Delta_0} Y_j^\mu(t\wedge
\sigma') > \zeta(1+\xi/4) \big|
Y^\mu_j(t_0\wedge\sigma') \le\zeta\xi/2 \Bigr) \le\frac{\xi
/2}{1+\xi/4} \leq\xi/2.
\end{eqnarray*}
Since $\zeta, \xi$ were arbitrary, this proves that $Y_j(t) \to
y_j(t)$ in probability uniformly on $[0,\Delta_\varepsilon\wedge
\sigma']$.
Since we have already shown that $P(\sigma_0 \vee\sigma_1 \leq
\sigma_2 ) \to0$, this completes the proof.
\end{pf*}

\subsection{Birth of a new type}

In this section, we prove Proposition \ref{case1}. Note that
$y_k(\Delta_\varepsilon) = 1-\varepsilon$ for small $\varepsilon$
since $\Delta_0=\delta
_k^0$ and choose $\bar{\varepsilon} = \bar{\varepsilon}(y^0)$ small
enough so that
the limiting dynamical system satisfies $y_j(\Delta_0+t) < \alpha+
t\rho/\gamma-2\eta$, $j \neq m,$ and $y_{j-1}(t) - y_j(t) < 1-2\eta$
for all $j \geq0$, $t \leq\bar{\varepsilon}$ and $\eta$ sufficiently
small. Since the result for $j \neq k$ follows from the arguments used
to prove Proposition \ref{interior}, we only need to prove \eqref
{kplus1}. To explain these inequalities we note that our limiting
process has $y_{k+1}(\Delta_0 +t) =\lambda_{k+1-m}t/\gamma$ and
$y_{k+1}(\Delta_0-t) = 0$ for small $t$. However, when $t$ is small
the number of type-($k+1$) individuals is small and deterministic
approximations are not valid. The best we can do is to say that
$Y_{k+1}^\mu(t)$ cannot get too far above the line with slope $\lambda
_{k+1-m}/\gamma$ that starts just before time $\Delta_0$ [the first
inequality in \eqref{kplus1}] or too far below the line with slope
$\lambda_{k+1-m}/\gamma$ that starts just after time $\Delta_0$ (the
second inequality).

We begin by defining branching processes $Z_{k+1,a}^\mu(t)$, $a=u,\ell
,$ with initial populations $Z_{k+1,a}^\mu(0)=0$ and per particle
birth and death rates given by Table \ref{table:compbp} in Section
\ref{sec-startup}, but with immigrations at rate $e^{\lambda_{k-m}t}
$. The methods used in the proof of the next result closely parallel
the methods used to prove Lemmas \ref{Zjmart}--\ref{Zmgcor} in
Section~\ref{sec-startup}.


\begin{lemma} \label{bploglim}
Let $0 < t_1 < t_2$. Then,
\[
Y_{k+1,a}^\mu(t) \equiv(1/L) \log^+ Z_{k+1,a}(Lt/\gamma) \to
t\lambda
_{k+1-m}/\gamma
\]
in probability uniformly on $[t_1,t_2]$ as $\mu\to0$ for $a=u,\ell$.
\end{lemma}

\begin{pf}
We prove the result for $a=u$ and drop the subscripts $u$ from all
quantities. For ease of notation, we will also write $Z(t) =
Z_{k+1}^\mu(t)$ but leave the $\mu$ superscript on $\lambda
_{k+1}^\mu$ to distinguish it from $\lambda_{k+1} = (1+\rho
)(1+\gamma
)^k-1$. Notice that
\[
\lambda_{k+1}^\mu= \frac{(1+\rho)(1+\gamma)^{j-m}}{1+\gamma_{-m} m
\mu^\eta} - \frac{1-(k+1) \mu^\eta}{1+\gamma_{k+1-m} (k+1-m) \mu
^\eta} \to\lambda_{k+1-m}
\]
as $\mu\to0$.

Define $M(t) = e^{-\lambda_{k+1}^\mu t}Z(t) - E(e^{-\lambda_{k+1}^\mu
t}Z(t))$. The same argument in the proof of Lemma \ref{Zjmart} then
implies that $M(t)$ is a martingale (with respect to the $\sigma
$-algebra generated by $Z(s)$, $s \leq t$). Furthermore, we have
%
\begin{equation}\label{meanbnd2}
E(e^{-\lambda_{k+1}^\mu t}Z(t)) = \bigl(1-e^{-(\lambda_{k+1}^\mu- \lambda
_{k-m})t}\bigr)/(\lambda_{k+1}^\mu- \lambda_{k-m})
\end{equation}
and a similar argument to the one used to prove \eqref{varZjt} in
Section \ref{startup} implies that
\begin{eqnarray*} \label{varbnd2}
\operatorname{var}(e^{-\lambda_{k+1}^\mu t}Z(t)) & \leq C.
\end{eqnarray*}
From the $L_2$ maximum inequality and Chebyshev's inequality, we can
conclude that
\begin{eqnarray*}
P\Bigl(\sup_{0\leq s \leq t_2} M(s) > L^{1/2} \Bigr) \to0
\end{eqnarray*}
as $\mu\to0$. This yields a result analogous to Lemma \ref{Zmgbnd}
in Section \ref{sec-startup}. The conclusion of Lemma \ref{bploglim}
then follows using the same argument as in the proof of Lemma~\ref
{Zmgcor} since $\lambda_{k+1}^\mu\to\lambda_{k+1-m} > \lambda
_{k-m}$ as $\mu\to0$ and so \eqref{meanbnd2} implies that there
exist $c_1,c_2 >0$ such that
\begin{eqnarray*}
c_1 \leq\sup_{t \leq s} E(e^{-\lambda_{k+1}^\mu s}Z(s)) \leq c_2
\end{eqnarray*}
for all $t >0$ if $\mu$ sufficiently small.
\end{pf}

\begin{pf*}{Proof of Proposition \ref{case1}} Suppose that $\varepsilon<
\varepsilon_1 =
\bar{\varepsilon} \wedge1/(2(\lambda_{k+1-m}/\lambda_{k-m}+\lambda
_{k+1-m}/\gamma))$.
Let $\eta>0$ and define
\begin{eqnarray*}
A_1 &=& \{|Y_k^\mu(t) - y_k(t)| \leq\varepsilon/4,   Y_j^\mu(t) <
\alpha
+\rho t -\eta,   \forall  j \leq k,   j \neq m,   t \leq\Delta
_0+\varepsilon\},
\\
A_2 & =& \{T_{k+1}^\mu\geq\Delta_{3\varepsilon/4}L/\gamma\},
\\
A_3 &=& \{  T_{k+2}^\mu\geq(\Delta_0+\varepsilon)L/\gamma,
Y_{k+1}^\mu(t)
< 1 -\eta,   \forall  t \leq\Delta_0+\varepsilon\},
\\
A &=& A_1 \cap A_2 \cap A_3.
\end{eqnarray*}
Note that $P(A_1^c) \to0$ by Propositions \ref{interior} and
 \ref{case1} applied to $j \neq k,$ while $P(A_2^c) \to0$
by Lemma \ref{nonewguys}. Recalling that
\[
y_k(\Delta_{3\varepsilon/4}+t) = 1-3\varepsilon/4 + t\lambda
_{k-m}/\gamma,
\]
we have
\begin{eqnarray*}
X_k^\mu\bigl(L(\Delta_{3\varepsilon/4}+t)/\gamma\bigr)e^{-\lambda_{k-m}
Lt/\gamma} \leq
1/\mu
\end{eqnarray*}
for all $t \leq(\Delta_0+\varepsilon)-\Delta_{3\varepsilon/4}$ on
$A$. Therefore,
using the bounds on the birth and death rates given in Lemma \ref
{bdbnds}, we can couple $X_{k+1}^\mu(L(\Delta_{3\varepsilon
/4}+t)/\gamma)$ with
$Z_{k+1,u}^\mu(Lt/\gamma)$, $a=u,$ in a similar manner to Lemma \ref
{coupling} so that on $A$, we have
\[
X_{k+1}^\mu\bigl(L(\Delta_{3\varepsilon/4}+t)/\gamma\bigr) \leq Z_{k+1,u}^\mu
(Lt/\gamma)
\]
for all $t \leq\Delta_0+\varepsilon-\Delta_{3\varepsilon/4}$.
Lemma \ref
{bploglim} then yields the first limit in \eqref{kplus1}, provided we
are on $A$. However, we then have
\[
P\bigl( \{Y_{k+1}^\mu(t) \leq1-\eta,   \forall  t \leq\Delta
_0+\varepsilon
\} \cap A\bigr) \to1,
\]
and it follows from this that $P(A_3^c)\to0$ for small $\eta$, which
proves that the first limit in \eqref{kplus1} holds. To prove the
second limit, we use the fact that
\[
y_k\bigl(\Delta_0 + (\varepsilon/4)(\gamma/\lambda_{k-m})+t\bigr) =
1+\varepsilon/4 + t\lambda
_{k-m}/\gamma
\]
to conclude that
\begin{eqnarray*}
X_k^\mu\bigl(L\bigl(\Delta_0 + (\varepsilon/4)(\gamma/\lambda
_{k-m})+t\bigr)/\gamma\bigr)e^{-\lambda
_{k-m} Lt/\gamma} \geq1/\mu
\end{eqnarray*}
for all $t \leq\varepsilon(1-\gamma/(4\lambda_{k-m}))$ on $A$.
Hence, we can
couple $X_{k+1}^\mu(L(\Delta_0 + (\varepsilon/4)(\gamma/\lambda
_{k-m})+t)/\gamma
)$ with $Z_{k+1,\ell}^\mu(Lt/\gamma)$ so that
\[
X_{k+1}^\mu\bigl(L\bigl(\Delta_0 + (\varepsilon/4)(\gamma/\lambda
_{k-m})+t\bigr)/\gamma\bigr) \geq
Z_{k+1,\ell}^\mu(Lt/\gamma)
\]
for all $t \leq\varepsilon(1-\gamma/(4\lambda_{k-m}))$ on $A$ and
the second
part of \eqref{kplus1} again follows from Lemma \ref{bploglim}.
\end{pf*}

\subsection{Change in dominant type}

In this section, we prove Proposition \ref{case2}. We begin with some
notation. Let
\begin{eqnarray*}
\hat{y}_j(t) \equiv
\cases{ (y_j^0 + \lambda_{j-m}t/\gamma)^+, &\quad $\mbox{if } j
\leq k$,
\cr
0, &\quad $\mbox{if } j > k.$
}
\end{eqnarray*}
Note that $y_j(t) \leq\hat{y}_j(t)$ for all $t$ with equality if $t
\leq\Delta_0$. As in Section 4.1, let $a_j(t) = \alpha+\rho t$ if $j
\neq k$, $a_k = 1$ and choose $\bar{\varepsilon} = \bar{\varepsilon
}(y^0)$ and
$\eta>0$ so that: (i) $\hat{y}_j(t) < a_j(t) - 2\eta$ for all $t
\leq\Delta_0+\bar{\varepsilon}$, $j \neq m,n;$ (ii) $y_{j-1}(t) -
y_j(t) <
1-2\eta$ for all $t \leq\Delta_0+\bar{\varepsilon}$; (iii) $y_j(t)
\geq
\alpha+\rho t - \eta/4$ for all $\Delta_{\bar{\varepsilon}} \leq t
\leq
\Delta_0 + \bar{\varepsilon}$, $j=m,n$. Without loss of generality, suppose
that $\eta< (\alpha+ \Delta_\varepsilon\rho/\gamma)/4$. Let
$\sigma_i(j)$,
$i=0,1,2,3,$ be as in Section 4.1 and set
\[
\bar{\sigma}(j) = \sigma_0(k+1) \wedge\Bigl(\min_{i \neq m,n} \sigma
_1(i)\Bigr) \wedge\sigma_2(j) \wedge\sigma_3(j).
\]
Our first lemma sets the stage for the battle between types $m$ and $n$
by showing that all other types remain smaller than these two.

\begin{lemma} \label{infmeanub}
If $j \neq m,n$ and $\xi>0,$ then
\[
P\Bigl(\sup_{t \leq\bar{\sigma}(j)} \bigl(Y_j^\mu(t) - \hat{y}_j(t)\bigr)
> \xi\Bigr) \to0.
\]
\end{lemma}

\begin{pf}
This follows directly from Lemma \ref{martbnd} and (I)--(IV) from
Lemma \ref{infmeanlim}.
\end{pf}

 Note that Lemma \ref{infmeanub}, Lemma \ref{poplim} and our choice
of $\bar{\varepsilon}$ imply that
%
\begin{eqnarray} \label{Xjsmall}
\sup_{t \leq\Delta_0+\bar{\varepsilon}} \frac{X_j(Lt/\gamma
)}{N^\mu(t)}
&\leq\mu^{\eta}
\end{eqnarray}
for all $j \neq n,m$ with high probability. Furthermore, $X_k(Lt/\gamma
) <
(1/\mu)^{1-\eta}$ for all $t \leq\Delta_0+\bar{\varepsilon}$ with high
probability, and hence Lemma \ref{nonewguys} implies that as $\mu\to0,$
%
\begin{equation}\label{nonewguys2}
P\bigl(T_{k+1}^\mu\leq(\Delta_0+\bar{\varepsilon}) L/\gamma\bigr) \to0.
\end{equation}

%
%

Let
\begin{eqnarray*}
R_j^\mu(t) \equiv\frac{X_j^\mu(L\Delta_\varepsilon/\gamma
+t)}{N^\mu(t)}
\end{eqnarray*}
be the fraction of $j$'s in the population at times greater than
$L\Delta_\varepsilon/\gamma$. Then, as a consequence of \eqref
{Xjsmall} and
\eqref{nonewguys2}, we have
%
\begin{eqnarray} \label{RnRmsum}
0 \leq1 - \bigl(R_m(Lt/\gamma) + R_n(Lt/\gamma)\bigr) \leq(k+1)\mu^{\eta}
\end{eqnarray}
for all $t \leq\Delta_0+\bar{\varepsilon}-\Delta_\varepsilon$ on
a set $A$ with
$P(A^c) \to0$. Our next result concerns the change of power from $m$'s
to $n$'s. To state the result, let
\[
f(r) \equiv r(1-r) \frac{\lambda_{n-m}}{1+\gamma_{n-m}r}
\]
and define $r_j^\mu(t)$, $j=n,m,$ as the (random) solutions to the
initial value problem
\begin{eqnarray} \label{ivpsol}
\frac{dr^\mu_m}{dt} &=& -f(1-r_m^\mu) \equiv f_m(r_m^\mu), \nonumber
\\[-8pt]\\[-8pt]
\frac{dr^\mu_n}{dt} &=& f(r_n^\mu) \equiv f_n(r_n^\mu),\nonumber
\end{eqnarray}
with initial conditions $r^\mu_j(0) = R_j^\mu(0)$, $j=n,m$.

\begin{lemma} \label{kurtz}
There exists an $\varepsilon_2 = \varepsilon_2(y^0) >0$ such that for $j=n,m$,
\begin{eqnarray*}
P\Bigl(\sup_{t \leq\Delta_0-\Delta_\varepsilon+\varepsilon}
|R_j^\mu
(Lt/\gamma)-r_j^\mu(Lt/\gamma)|> \mu^{\eta/2} \Bigr) \to0
\end{eqnarray*}
as $\mu\to0$ for all $\varepsilon< \varepsilon_2$.
\end{lemma}

\begin{pf}
We will prove the result by calculating the infinitesimal mean and
variance of $R_j^\mu(t)$. Without loss of generality, we assume that
$\varepsilon< \bar{\varepsilon}$ so that by \eqref{Xjsmall}--\eqref{RnRmsum}, we have
\begin{eqnarray} \label{Rcond}
\sum_{j \neq n,m} R_j^\mu(t) &\leq&(k+1) \mu^\eta, \nonumber
\\[-8pt]\\[-8pt]
0 &\leq&1-\bigl(R_n^\mu(t) + R_m^\mu(t)\bigr)  \leq(k+1)\mu^\eta\nonumber
\end{eqnarray}
for all $t \leq L(\Delta_0 + \varepsilon- \Delta_{\varepsilon
})/\gamma$ on a set $A$
with $P(A^c) \to0$ as $\mu\to0$. Note also that Lemma \ref{poplim},
the fact that $N^\mu$ is nondecreasing and our choice of $\eta<
(\alpha+\Delta_\varepsilon\rho/\gamma)/4$ together imply that
%
\begin{equation}\label{popsize}
N^\mu(t) \geq C (1/\mu)^{\alpha+\Delta_\varepsilon\rho/\gamma
-\eta} \geq C
(1/\mu)^{3\eta} \qquad\forall  t \geq L\Delta_\varepsilon/\gamma
\end{equation}
on a set $A$ with $P(A^c) \to0$ as $\mu\to0$. We will therefore
assume that the inequalities in \eqref{Rcond} and \eqref{popsize}
hold for the remainder of the proof and write $O(\mu^\eta)$ for any
quantity whose absolute value is bounded above by $C \mu^\eta$
uniformly for $t \leq L(\Delta_0 + \varepsilon- \Delta_{\varepsilon
})/\gamma$ on a set
$A$ with $P(A^c) \to0$ as $\mu\to0$. It all also convenient to write
\[
c_\varepsilon= \Delta_0 + \varepsilon- \Delta_{\varepsilon} =
\bigl(1+\gamma(\lambda_{n-m}-\rho
)^{-1}\bigr)\varepsilon.
\]

By looking at the rates for the chain $(N^\mu(t),X^\mu(t))$, the
fraction $R_j^\mu(t)$ has the following jump rates corresponding to
the events $x_j/N \mapsto(x_j+1)/N$, $x_j/N \mapsto(x_j-1)/N$, $x_j/N
\mapsto(x_j+1)/(N+1)$ and $x_j/N \mapsto x_j/(N+1)$, respectively:
\begin{eqnarray*}
r_j &\mapsto& r_j +1/N \qquad\mbox{rate: }  N(1 - r_j) \frac
{(1+\gamma)^j r_j}{w} + \mu N r_{j-1};
\\
r_j &\mapsto& r_j -1/N \qquad\mbox{rate: }  Nr_j \frac
{w-(1+\gamma
)^j r_j}{w } + \mu N r_j;
\\
r_j &\mapsto& r_j + (1-r_j)/(N+1) \qquad\mbox{rate: } \rho N
\frac{ (1+\gamma)^j r_j}{w };
\\
r_j &\mapsto& r_j - r_j/(N+1) \qquad\mbox{rate: } \rho N \frac
{w - (1+\gamma)^j r_j}{w },
\end{eqnarray*}
where $w \equiv\sum_{i \geq0} (1+\gamma)^i r_i$. From these expressions
for the rates, we can see that the infinitesimal mean of $R_j^\mu$ is
given for $r \in\mathcal{S}^N/N$ by
\begin{eqnarray*}
B_j(r) &=& \frac{(1+\gamma)^j r_j(1 - r_j)}{w} + \mu r_{j-1} - \biggl(
\frac{r_j(w-(1+\gamma)^j r_j)}{w } + \mu r_j\biggr)
 \\
&&{} + \frac{\rho N}{N+1} \frac{ (1+\gamma)^j r_j(1-r_j)}{w } -
\frac
{\rho N}{N+1} \frac{r_j(w - (1+\gamma)^j r_j)}{w }
\\
& =& \biggl(1+\frac{\rho N}{N+1}\biggr)\biggl( \frac{r_j((1+\gamma
)^j -
w)}{w} \biggr) + \mu(r_{j-1}-r_j).
\end{eqnarray*}
Similarly, the infinitesimal variance is given by
\begin{eqnarray*}
A_j(r) &=& \frac{1}{N} \biggl(\frac{(1+\gamma)^j r_j(1 - r_j)}{w} +
\mu
r_{j-1} + \frac{r_j(w-(1+\gamma)^j r_j)}{w } + \mu r_j \biggr)
\\
&&{} + \frac{\rho N}{(N+1)^2} \biggl(\frac{ (1+\gamma)^j
r_j(1-r_j)^2}{w } + \frac{r_j^2(w - (1+\gamma)^j r_j)}{w }\biggr)
 \\
&=& \frac{1}{N}\biggl( \biggl(1+ \frac{\rho N^2}{(N+1)^2}\biggr)
\frac{r_j((1+\gamma)^j-2(1+\gamma)^j r_j+w)}{w}
 \\
&&{}\hspace*{88pt}  - \frac{\rho r_j(1-r_j) N^2}{(N+1)^2} + \mu
(r_{j-1}+r_j)\biggr),
\end{eqnarray*}
where, in the second line, we have added and subtracted $\rho r_j
N/(N+1)^2$ from the first. Note that \eqref{popsize} and the fact that
$r_j \in[0,1]$ together imply
that
%
\begin{equation}\label{varmu}
A_j^\mu(R^\mu(s)) = O(\mu^{3\eta})
\end{equation}
for all $s \leq Lc_\varepsilon/\gamma$.

Now, \eqref{Rcond} implies that
\begin{eqnarray*}
w(R^\mu(s)) &=& (1+\gamma)^m R_m^\mu(s) + (1+\gamma)^n R_n^\mu(s) +O(\mu^\eta)
\\
& =& (1+\gamma)^m[1+ \gamma_{n-m}r_n^\mu(s)] + O(\mu^\eta)
\end{eqnarray*}
for all $s \leq L c_\varepsilon/\gamma$ and hence
\begin{eqnarray*}
&&B_n(R^\mu(s))
\\
&&\qquad = \biggl(1+\frac{\rho N^\mu(L\Delta_\varepsilon/\gamma+s)}{N^\mu
(L\Delta
_\varepsilon/\gamma+s)+1}\biggr)
\\
&&{}\quad\qquad\times R_n^\mu(s)\biggl( \frac{(1+\gamma
)^n - (1+\gamma
)^m[1+ \gamma_{n-m}R_n^\mu(s)] + O(\mu^\eta) }{(1+\gamma)^m[1+
\gamma
_{n-m}R_n^\mu(s)] + O(\mu^\eta) } \biggr)
+ O(\mu^\eta)
\\
&&\qquad = \biggl(1+\frac{\rho N^\mu(L\Delta_\varepsilon/\gamma+s)}{N^\mu
(L\Delta
_\varepsilon/\gamma+s)+1}\biggr) R_n^\mu(s)\biggl( \frac{\gamma
_{n-m}-\gamma
_{n-m}R_n^\mu(s)}{1+\gamma_{n-m}R_n^\mu(s)} \biggr) + O(\mu^\eta)
\\
&&\qquad = f_n(R_n^\mu(s)) + O(\mu^\eta)
\end{eqnarray*}
for all $s \leq L c_\varepsilon/\gamma$, the last equality following
from \eqref
{popsize} and the definition of $f_n$. Similarly, writing
\begin{eqnarray} \label{wexp}
w(R^\mu) &=& (1+\gamma)^m R_m^\mu+ (1+\gamma)^n R_n^\mu+ O(\mu
^\eta)\nonumber
\\[-8pt]\\[-8pt]
& =& (1+\gamma)^m[1+ \gamma_{n-m}(1-R_m^\mu)] + O(\mu^\eta),\nonumber
\end{eqnarray}
we obtain
\begin{eqnarray*}
B_m(R^\mu(s)) = f_m(R_m^\mu(s)) + O(\mu^\eta)
\end{eqnarray*}
for all $s \leq Lc_\varepsilon/\gamma$. Combining this with \eqref
{varmu}, the
fact that $|f'(r_j)| \leq\gamma(1+\gamma)$ for all $r_j \in[0,1]$
and the
proof of Theorem 2.11 in \cite{TK}, we obtain the result. (Theorem
2.11 in \cite{TK} applies directly if we replace $r_j^\mu(t)$ with
$r_j(t)$, the solution to \eqref{ivpsol} with initial conditions
$r_j(0) = \lim_{\mu\to0} R_j^\mu(0)$, but it is easy to see that
the same proof applies if we use $r_j^\mu(t)$ since $r_j^\mu(t)$ is
the solution to \eqref{ivpsol} with random initial conditions $R_j^\mu(0)$.)
\end{pf}


The next step is to analyze the differential equations for $j=m,n$ in
Lemma \ref{kurtz}. We will carry out the analysis for $j=n$ [for
$j=m$, apply the analysis below to $1-r_m^\mu(s)$]. To begin, write
\[
r_n^\mu(t) = \frac{X_n^\mu(L\Delta_\varepsilon/\gamma)}{N^\mu
(L\Delta_\varepsilon
/\gamma)}\exp\biggl\{\int_0^t g_n(r_n^\mu(s)) \,ds\biggr\}
\]
with
\[
g_n(r) \equiv\frac{\lambda_{n-m}(1-r) }{1+\gamma_{n-m} r}.
\]
Note that we have the following set of bounds on the growth rate $g_n$:
%
\begin{eqnarray} \label{rnbnds}
(1 - L^{-2}) \frac{\lambda_{n-m}}{ 1 + \gamma
_{n-m}L^{-2}} & \le& g_n(r_n^\mu) \le \lambda_{n-m} \qquad \mbox{when
} r_n^\mu< L^{-2};\nonumber
\\
(1 - r^\mu_n) \frac{\lambda_{n-m} }{ 1 + \gamma_{n-m}}& \le&
g_n(r_n^\mu) \le(1-r^\mu_n) \lambda_{n-m}\nonumber
\\[-8pt]\\[-8pt]
&&\eqntext{\qquad \mbox{when } L^{-2}
\le r_n^\mu\le1 - L^{-2};}
\\
0 & \le& g_n(r_n^\mu) \le L^{-2} \lambda_{n-m} \qquad \mbox{when }
r_n^\mu\ge1- L^{-2}.\nonumber
\end{eqnarray}

\begin{lemma} \label{quickgrow}
Let $s^\mu_1 = \inf\{ s \dvtx r_n^\mu(s) \ge L^{-2} \}$ and $s^\mu_2 =
\inf\{ s\dvtx r_n^\mu(s) \ge1- L^{-2} \}.$
We then have $s^\mu_i/(L/\gamma) \to\Delta_0$ for $i=1,2$ and
$(s_2^\mu- s_1^\mu)/L \to0$ as $\mu\to0$.
\end{lemma}

\begin{pf}
Let
\[
r_\ell^\mu(s) = \frac{X_n^\mu(L\Delta_\varepsilon/\gamma)}{N^\mu
(L\Delta
_\varepsilon/\gamma)} e^{\lambda_{n-m} c_\mu s} \quad\mbox{and}
\quad
r_u^\mu(s) = \frac{X_n^\mu(L\Delta_\varepsilon/\gamma)}{N^\mu
(L\Delta_\varepsilon
/\gamma)} e^{\lambda_{n-m} s},
\]
where $c_\mu= (1 - L^{-2})(1 + \gamma_{n-m}L^{-2})^{-1}$. It is then
clear from the first bound in \eqref{rnbnds} that
\[
r_\ell^\mu(s) \leq r_n^\mu(s) \leq r_u^\mu(s)
\]
for all $s \leq s_1^\mu$. Since $Y_n(\Delta_\varepsilon) \to
y_n(\Delta_\varepsilon
)$, $F^\mu(\Delta_\varepsilon) \to\alpha+\Delta_{\varepsilon
}\rho/\gamma$ by
Proposition \ref{interior}, letting $s_\ell^\mu$ and $s_u^\mu$ be
the times that $r_\ell^\mu$ and $r_u^\mu$ hit $L^{-2}$, we have
$s_a^\mu/(L/\gamma) \to\Delta_0$ as $\mu\to0$ for $a=\ell,u,$ which
proves the result for $i=1$. To prove the result for $i=2$ we use the
bounds in the second line of \eqref{rnbnds} along with the fact that
the logistic $dx/dt = \beta x(1-x)$ rises from $L^{-2}$ to $1-L^{-2}$
in time $(4/\beta) \log L,$ to conclude that
\[
\frac{s_2^\mu- s_1^\mu}{L} \leq\frac{C \log L}{ L} \to0
\]
as $\mu\to0,$ which completes the proof.
\end{pf}

\begin{lemma} \label{rjlim}
\[
(1/L) \log^+ [N^\mu(Lt/\gamma) r_j^\mu(Lt/\gamma)] \to y_j(\Delta
_\varepsilon+ t)
\]
uniformly on $[0,T]$ for any $T >0$, $j=n,m$.
\end{lemma}

\begin{pf}
We prove the result for $j=n$. Write
\begin{eqnarray*}
&&(1/L) \log^+ [N^\mu(Lt/\gamma) r_j^\mu(Lt/\gamma)] - y_j(\Delta
_\varepsilon+t)\nonumber
\\
&&\qquad= [Y_j^\mu(\Delta_\varepsilon) - y_j(\Delta_\varepsilon)]
\\
&&{}\qquad \quad+ (1/L) \biggl[ \int_0^{Lt/\gamma} \bigl(g_n(r_n^\mu(s)) - \ell
_n(s)\bigr) \,ds\biggr],\nonumber
\end{eqnarray*}
where $\ell_n(s) = \lambda_{n-m} \mathbf{1}_{s \leq(\Delta
_0-\Delta_\varepsilon)}$. The first term in brackets converges to 0 in
probability by Proposition \ref{interior}. To control the second term,
split up the integral as
\begin{eqnarray*}
\int_0^{tL/\gamma} = \int_0^{s_1^\mu\wedge t} + \int_{s_1^\mu
\wedge
t}^{s_2^\mu\wedge t} + \int_{s_2^\mu\wedge t}^t.
\end{eqnarray*}
Using the bounds in \eqref{rnbnds} and applying Lemma \ref
{quickgrow}, we conclude that each of these integrals is $o(L),$ which
yields the result.
\end{pf}

\begin{pf*}{Proof of Proposition \ref{case2}} Let $\varepsilon<
\varepsilon_2$ and
suppose first that $j=m,n$. Writing
\begin{eqnarray*}
Y_j^\mu(t) - y_j(t)
& =& (1/L)\log^+ R_j^\mu(Lt/\gamma) - (1/L)\log^+r_j^\mu(Lt/\gamma
)
\\
&&{}+
(1/L) \log^+ [N^\mu(Lt/\gamma) r_j^\mu(Lt/\gamma)] - y_n(t),
\end{eqnarray*}
we can see that since
\[
r_j^\mu(0) \geq\frac{X_n^\mu(L\Delta_\varepsilon/\gamma)}{N^\mu
(L\Delta
_\varepsilon/\gamma)} \geq\mu^{\eta/3}
\]
for all $t \geq0$ with high probability as $\mu\to0$ by our choice
of $\bar{\varepsilon}$, the result follows from Lemmas \ref{kurtz}
and \ref
{rjlim}. Suppose now that $j \neq m,n$. If $j > k$, the result follows
from \eqref{nonewguys2}, so it remains to prove the result for $j \leq
k$. In view of Lemma \ref{martbnd}, it suffices to show that
\[
P\biggl(\sup_{t \leq(\Delta_0+\varepsilon) \wedge\bar{\sigma
}(j)}
\bigg|y_j^0 + \int_0^t B_j(Y_j^\mu(s))\,ds - y_j(t) \bigg| > \xi\biggr)
\to0
\]
and then follow the argument from the proof of Proposition \ref
{interior} to yield the result. However, now that we have proven that
Proposition \ref{case2} holds for $j=m,n$, we can argue as in the
proof of Lemma \ref{infmeanlim} to conclude that
\begin{eqnarray*}
g_{j,1}(Y^\mu(t)) \to
\cases{ \dfrac{(1+\rho)(1+\gamma)^{j-m}}{\gamma}, &\quad$ 0\le t <
\Delta
_0,$
 \cr
\dfrac{(1+\rho)(1+\gamma)^{j-n}}{\gamma}, &\quad$ \Delta_0 < t \le
(\Delta_0+\varepsilon
)$,
}
\end{eqnarray*}
uniformly on compact subsets of $[0,(\Delta_0+\varepsilon) \wedge
\bar
{\sigma}(j)] - \{\Delta_0\}$. This replaces the second part of (IV)
from the proof of Lemma \ref{infmeanlim}, and the result follows after
using (I)--(III) from the proof of Lemma \ref{infmeanlim}.
\end{pf*}


\printaddresses

\end{document}